\documentclass[leqno,12pt]{amsart}
\usepackage[letterpaper,margin=1in]{geometry}
\usepackage{amssymb}
\usepackage[usenames, dvipsnames]{color}
\usepackage{enumitem}

\usepackage{mathrsfs}

\usepackage{marginnote}
\usepackage{todonotes}
\usepackage{hyperref}
\usepackage{esint}
\usepackage{verbatim}
\usepackage[normalem]{ulem}

\theoremstyle{plain}
\newtheorem{theorem}{Theorem}[section]
\newtheorem{lemma}[theorem]{Lemma}

\newtheorem{proposition}[theorem]{Proposition}

\theoremstyle{definition}
\newtheorem{definition}[theorem]{Definition}

\theoremstyle{remark}
\newtheorem{remark}[theorem]{Remark}

\DeclareMathOperator*{\diver}{div}
\DeclareMathOperator{\bA}{\mathbf{A}}
\DeclareMathOperator*{\bB}{\mathbf{B}}

\DeclareMathOperator*{\bF}{\mathbf{F}}
\DeclareMathOperator*{\bG}{\mathbf{G}}
\DeclareMathOperator*{\F}{F}
\DeclareMathOperator*{\A}{A}

\DeclareMathOperator*{\bK}{\mathbb{K}}

\numberwithin{equation}{section}

\newcommand{\bN}{\mathbb{N}}

\newcommand{\bR}{\mathbb{R}}

\newcommand\p{\partial}

\def\Xint#1{\mathchoice
	{\XXint\displaystyle\textstyle{#1}}%
	{\XXint\textstyle\scriptstyle{#1}}%
	{\XXint\scriptstyle\scriptscriptstyle{#1}}%
	{\XXint\scriptscriptstyle\scriptscriptstyle{#1}}%
	\!\int}
\def\XXint#1#2#3{{\setbox0=\hbox{$#1{#2#3}{\int}$ }
		\vcenter{\hbox{$#2#3$ }}\kern-.6\wd0}}

\def\fint{\Xint-}

\makeatother

\begin{document}

\title[ Singular-Degenerate Equations with Solution-Dependent Coefficients]{Gradient Estimates for Singular-Degenerate Elliptic Equations with Solution-Dependent Coefficients}

\author[V. Le]{Van Le}
\address[V. Le]{Department of Mathematics, University of Tennessee, 227 Ayres Hall,
	1403 Circle Drive, Knoxville, TN 37996-1320, USA}
\email{nle12@vols.utk.edu}

\author[T. Phan]{Tuoc Phan}
\address[T. Phan]{Department of Mathematics, University of Tennessee, 227 Ayres Hall,
	1403 Circle Drive, Knoxville, TN 37996-1320, USA}
\email{phan@utk.edu}


\begin{abstract} 
We study a conormal boundary value problem for a class of singular-degenerate quasilinear elliptic equations in divergence form in bounded domains. The leading terms are measurable in the spatial variable and depend nonlinearly on both the solution and its gradient. The singularity or degeneracy is governed by a power of the distance to the domain boundary, with exponent greater than negative one. In particular, the corresponding weights need not belong to the Muckenhoupt $A_2$ class. Assuming that the weighted mean oscillations of the leading terms in the spatial variable are sufficiently small on small balls and that the domains have Lipschitz boundaries with sufficiently small Lipschitz constants, we establish global weighted $W^{1,p}$-estimates for weak solutions with bounded weighted mean oscillation. Our results extend recent estimates for equations whose leading terms are independent of the solution itself.
\end{abstract}


\maketitle
\vspace{-0.7cm}
\section{Introduction}
\subsection{Settings}
Let $\Omega$ be a nonempty bounded connected open set in $\bR^n$ with Lipschitz boundary $\partial\Omega$ and $n\in \mathbb{N}$.
We investigate the following class of singular-degenerate quasilinear equations  subject to conormal boundary conditions:
\begin{equation}\label{main-eqs}
\left \{ 
\begin{array}{cccl}
		\diver(\mu(x)\bA(x,u,\nabla u)) &= &\diver(\mu(x)\bF(x)) & x \ \in \ \Omega,\\
	\mu(x)\big(\bA(x,u (x),\nabla u(x))-\bF(x)\big)\cdot\overrightarrow{\nu}(x) &= &0 & x \ \in \ \partial\Omega.
\end{array} \right.
\end{equation}
Here, $u: \Omega \rightarrow \bR$ is an unknown solution, $\bF: \Omega \rightarrow \bR^n$ is a given measurable vector field, $\overrightarrow{\nu}$ denotes the outward unit normal to $\partial\Omega$. Moreover, $\bA: \Omega \times \overline{\bK} \times \bR^n\rightarrow \bR^n$ is a given vector field, where $\bK \subset \bR$ is an open interval which is possibly equal to $\bR$.  Also, $\mu$ in \eqref{main-eqs} is positive and continuous in $\Omega$ and satisfies 
\begin{equation}\label{mu-weight}
	\mu(x)= 	\left\{\begin{array}{cl}
		\text{dist}(x,\p \Omega)^\alpha &\text{if }	\text{dist}(x,\p \Omega)\leq r_0,\\
		 1 &\text{if }	\text{dist}(x,\p \Omega)>2r_0
	\end{array}\right.
\end{equation}
where $r_0\in(0,\operatorname{diam}(\Omega)/2)$ is sufficiently small and $\alpha\in(-1,\infty)$.  The  conormal boundary condition in \eqref{main-eqs} is understood in the weak sense as in Definition \ref{weak-sol-def} below.

\smallskip
We assume that $\bA= \bA(x,z,\xi)$ is measurable with respect to $x\in \Omega$ for every $(z,\xi)\in \bK\times(\bR^n\setminus\{0\})$, Lipschitz continuous with respect to $z \in \bK$ for a.e. $x\in \Omega$ and every $\xi\in \bR^n\setminus\{0\}$, and differentiable with respect to $\xi$ in $\bR^n\setminus\{0\}$ for each $z\in \bK$ and for a.e. $x\in \Omega$. In addition, we assume that there exists a constant $\Lambda\geq 1$ such that  $\bA: \Omega \times \overline{\bK} \times \bR^n\rightarrow \bR^n$ satisfies the following ellipticity and growth conditions:
\begin{equation}\label{bA-growth}
\left\{	\begin{array}{ll}
|\xi-\zeta|^2 \leq \Lambda \langle \bA(x,z,\xi)-\bA(x,z,\zeta), \xi -\zeta \rangle & \quad \forall \  z\in \overline{\bK},  \ \xi,\zeta\in\bR^n \\
|\bA(x,z_1,\xi)-\bA(x,z_2,\xi)|\leq \Lambda  |z_1-z_2||\xi| &\quad \forall \  z_1,z_2\in \overline{\bK}, \ \xi \in\bR^n  \\
\bA(x,z,0)=0,\ |D_\xi \bA(x,z,\xi)|\leq\Lambda
        & \quad \forall \  z\in\overline{\bK}, \ \xi \in\bR^n \setminus\{0\}
        \end{array} \right.
\end{equation}
for a.e.  $x\in \Omega$. 

\smallskip
Under assumptions \eqref{mu-weight} and \eqref{bA-growth}, equation \eqref{main-eqs} displays two distinctive features. First, its leading term, $\mu(x) \bA(x, u(x), \nabla u(x))$, depends nonlinearly on both $u$ and $\nabla u$. Thus, it includes very large class of quasilinear elliptic equations, see \cite[Chapter 4]{Ladyzhenskaya-Uraltseva}. Second, the leading term is singular near $\partial\Omega$ when $\alpha < 0$ and degenerate when $\alpha > 0$. Crucially, since $\alpha > -1$ can be arbitrarily large, $\mu$ needs not belong to the Muckenhoupt $A_2$ class as commonly studied in the literature (e.g., \cite{CMP, Fabes, FKS, MP}). Finally, when $\alpha = 0$ and $\mu \equiv 1$, the equation reduces to a standard uniformly elliptic quasilinear equation in which $\bA$ depends nonlinearly on both $u$ and $\nabla u$.

\smallskip
This work establishes global regularity estimates in weighted Sobolev space $W^{1,p}(\Omega, \mu)$ for weak solutions of \eqref{main-eqs}. Theorem~\ref{main-thrm} shows that if $\bF\in L^p(\Omega,\mu)^n$ for some $p\in[2,\infty)$, the weighted mean oscillations of $\bA$ are sufficiently small, and $\partial\Omega$ has a sufficiently small Lipschitz constant, then
\begin{equation}\label{ineqs}
\Big(\int_{\Omega} |\nabla u(x)|^pd\mu(x)\Big)^{1/p} \leq N\Big(\int_{\Omega}|\bF(x)|^p d\mu(x)\Big)^{1/p}
\end{equation}
for any weak solution $u$ of \eqref{main-eqs} satisfying $[[u]]_{\textup{BMO}_{R_0}(\Omega,\mu)}\leq M$, where $N>0$ is independent of $u$ and $\bF$. The theorem extends the result of \cite{Dong-Phan-Sire}, in which $\bA = \bA(x, \nabla u)$ is independent of $u$. When $\mu\equiv1$, it also extends and recovers unweighted $W^{1,p}$ estimates for quasilinear elliptic equations; see, for example, \cite{BPS, NP, P}.
\subsection{Motivation and related literature} Equation \eqref{main-eqs} describes a stationary diffusion problem with flux
\[
\mathbf{J}=\mu(x)\big(\bF(x)-\bA(x,u,\nabla u)\big),
\qquad \diver\mathbf{J}=0\ \text{in }\Omega,
\qquad \mathbf{J}\cdot\overrightarrow\nu=0\ \text{on }\partial\Omega.
\]
The conormal boundary condition prescribes zero normal flux through the boundary. The factor $\mu(x)$ describes the spatial degeneracy or singularity of the diffusion, while the dependence of $\bA$ on $u$ allows the diffusion law to vary with the state. For example, $\bA(x,z,\xi)=a(x,z)\xi$ satisfies \eqref{bA-growth} if
\[
\Lambda^{-1}\leq a(x,z)\leq\Lambda,
\qquad |a(x,z_1)-a(x,z_2)|\leq\Lambda|z_1-z_2|.
\]
In this case, the diffusion coefficient is $\mu(x)a(x,u)$. Near $\partial\Omega$, it tends to zero when $\alpha>0$ and becomes unbounded when $\alpha<0$. Thus, even this simple example requires estimates that retain the boundary weight.

\smallskip
We present a few models in applied areas that motivate this study. Moosaie  \cite[(2.2)-(2.9)]{Moosaie} considers steady heat conduction in a cylindrical shell with outer radius $R_o$. The temperature $T$ satisfies
\[
\diver\big(k(r,T)\nabla T\big)=0, \quad \text{where} \quad
k(r,T)=\left(\frac{r}{R_o}\right)^m(\lambda_0-\lambda_1 T),
\]
where $r$ is the distance to the axis. The factor $(r/R_o)^m$ describes the radial variation of the material, and $\lambda_0-\lambda_1 T$ gives the temperature dependence of the conductivity. For the lake equations with varying depth $b$, see \cite[Section 2]{Camassa-Holm-Levermore}. When $b$ vanishes like $\textup{dist}(x,\partial\Omega)^a$ at the shore, with $0<a<1$, the coefficient $b^{-1}$ in the stream-function operator $\diver(b^{-1}\nabla\psi)$ is comparable near the shore to the weight in \eqref{mu-weight} with $\alpha=-a$, see \cite[(1.4), (2.4)]{Bresch-Metivier}.
In a porous lithium-ion electrode, the electrolyte diffusive flux is $-B(x)D(c)\partial_x c$, where $c$ is the electrolyte concentration and $B(x)$ accounts for the pore geometry, see \cite[(1)]{Richardson-battery}. For a spherical electrode particle, the lithium concentration satisfies $r^2\partial_t c=\partial_r(r^2D(c)\partial_r c)$, where $r$ is the distance to the particle center, see \cite[(1)]{Renganathan-White}.

\smallskip
In pure mathematics, a concrete source of such weighted operators is the extension problem for the fractional Laplacian. For $0<s<1$, Caffarelli and Silvestre \cite{CS-extension} represent $(-\Delta)^s v$ by
\[
\diver(y^{1-2s}\nabla U)=0\quad\text{in }\bR^{n-1}\times(0,\infty),
\qquad U(x',0)=v(x'),
\]
\[
-\lim_{y\downarrow0}y^{1-2s}\partial_yU(x',y)
=c_s(-\Delta)^s v(x'),\qquad c_s>0.
\]
The interior equation is the flat model considered here with $\alpha=1-2s$, $\bA(x,z,\xi)=\xi$, and $\bF=0$. This representation replaces a nonlocal equation by a local equation with a power weight and a conormal boundary condition.  The same operator occurs in models of phase transitions and dislocations in crystals. Cabr\'e and Sire \cite{Cabre-Sire-app} study $(-\Delta)^s v=f(v)$ through the extension equation and the nonlinear boundary condition
\[
-\lim_{y\downarrow0}y^{1-2s}\partial_yU(x',y)
=c_sf(U(x',0)).
\]
Caffarelli, Roquejoffre, and Sire \cite{CRS-free-boundary} use the weighted extension to study variational problems with free boundaries. These examples motivate regularity estimates for the weighted interior operator. Their boundary conditions are determined by the particular model and need not be the homogeneous conormal condition in \eqref{main-eqs}.

\smallskip
For uniformly elliptic equations, the Calder\'on-Zygmund theory concerns whether the integrability of the data implies the corresponding integrability of the gradient. Caffarelli and Peral \cite{CP} developed a perturbative approach to $W^{1,p}$ estimates for elliptic equations in divergence form. Byun and Wang \cite{Buyn-Wang2,Buyn-Wang} established global estimates for equations with small BMO coefficients or nonlinearities on Reifenberg flat domains. They also studied the conormal derivative problem under similar assumptions on the coefficients and the boundary \cite{Byun-Wang-Conormal}. For nonlinear elliptic Dirichlet problems, Mengesha and Phuc \cite{MN} obtained global weighted gradient estimates on Reifenberg flat domains. The common feature of these results is the use of comparison estimates and maximal-function techniques to obtain higher integrability of the gradient, without requiring continuity of the coefficients in the spatial variable. Note that these mentioned papers do not cover the considered equations \eqref{main-eqs} as the main terms in these works are independent on $u$.

\smallskip
Another approach to $L^p$ estimates for elliptic and parabolic equations with discontinuous coefficients was developed by Krylov \cite{Krylov, Krylov2007}, using the Fefferman-Stein sharp function theorem \cite{FeffermanStein}. The method starts with mean oscillation estimates for derivatives of solutions to equations with simpler coefficients. Perturbation then gives bounds for the sharp functions of these derivatives in terms of Hardy-Littlewood maximal functions of the data and the derivatives themselves. The Fefferman-Stein theorem and maximal-function estimates yield the desired $L^p$ bounds when the coefficient oscillations are sufficiently small. Kim and Krylov \cite{MR2338417} studied elliptic equations with coefficients measurable in one variable and VMO in the remaining variables. Dong and Kim \cite{DK11b, DK11} developed the theory for elliptic and parabolic systems with variably partially BMO coefficients and subsequently obtained estimates with $A_p$ weights \cite{MR3812104}.

\smallskip
Boundary value problems for degenerate elliptic operators were studied by Murthy and Stampacchia \cite{MuSt1}. Local regularity of solutions was investigated by Fabes \cite{Fabes}, and Fabes, Kenig, and Serapioni \cite{FKS} established local H\"older estimates and Harnack inequalities for a class of degenerate elliptic equations with $A_2$ weights. H\"older regularity for certain equations with boundary degeneracy was also studied in \cite{Wang2}. For gradient estimates, Cao, Mengesha, and Phan \cite{CMP} proved global weighted $W^{1,p}$ estimates for linear singular-degenerate elliptic Dirichlet problems with $A_2$ Muckenhoupt weights, while Mengesha and Phan \cite{MP} considered coefficients that singular-degenerate as an $A_2$ weight of one spatial variable. Related estimates for singular-degenerate elliptic and parabolic equations were developed in \cite{Dong-Phan,Dong-Phan-1,Dong-Phan-3}. More recently, Chang and Zheng \cite{CZ} established global estimates in weighted generalized Morrey spaces for linear singular-degenerate elliptic systems with Dirichlet boundary conditions. The work most directly related to the present paper is Dong, Phan, and Sire \cite{Dong-Phan-Sire}. They studied conormal problems with leading terms of the form $\mu(x)\bA(x,\nabla u)$, where $\mu$ behaves like $\operatorname{dist}(x,\partial\Omega)^\alpha$, and established weighted $W^{1,p}$ estimates for $\alpha\in(-1,\infty)$ under a small weighted mean oscillation assumption. In particular, their results allow weights outside the $A_2$ class.

\smallskip
Dependence of the coefficients on the solution introduces another issue. On one hand, quasilinear elliptic equations with coefficients depending on both the solution and its gradient are studied in the classical monograph of Ladyzhenskaya and Ural'tseva \cite[Chapter 4]{Ladyzhenskaya-Uraltseva} in which H\"{o}lder estimates are developed for bounded weak solutions. On the other hand, regularity theory in Sobolev spaces for this class of equation has just developed recently. For equations of the form $\diver\bA(x,u,\nabla u)=\diver\bF$, Nguyen and Phan \cite{NP} introduced a scaling parameter and established  interior $L^p$ gradient estimates for locally bounded weak solutions, with measurable dependence on $x$ and Lipschitz dependence on $u$. Phan \cite{P} obtained interior, boundary, and global estimates under a BMO assumption on the solution, in place of boundedness, for quasilinear equations with inhomogeneous conormal data. Weighted estimates for solution-dependent degenerate equations with $A_2$ degeneracy and Dirichlet boundary data were established in \cite{P2}. These results distinguish the oscillation of the coefficients in $x$ from the oscillation introduced by the solution. See also \cite{GN, BPS} for a few following up works in this research line which covers $p$-Laplacian type equations.

\smallskip
In this paper, we combine boundary singularity and degeneracy with dependence of $\bA$ on both $u$ and $\nabla u$. We assume a bound $M$ for the weighted BMO seminorm of $u$ and small weighted mean oscillation of $\bA$ in $x$. The bound $M$ need not be small; the required smallness of the coefficient oscillation  in $x$ depends on $M$. Under these assumptions, we prove
\[
\|\nabla u\|_{L^p(\Omega,\mu)}
\leq C\|\bF\|_{L^p(\Omega,\mu)},\qquad p\in[2,\infty).
\]
The result extends the solution-independent conormal theory in \cite{Dong-Phan-Sire} to solution-dependent coefficients under the hypotheses stated below.  The main result of the paper is stated in Theorem \ref{main-thrm} below. Corresponding local boundary regularity estimates and local interior regularity estimates are also established in Theorem \ref{flat-thrm} and Theorem \ref{interior-thrm}.
\subsection{Definitions and main result} Throughout the paper,  for $p\in[1,\infty)$ and for an open set $B\subset\bR^n$ and for a non-negative locally integrable function $\omega$ defined on $B$, $L^p(B, \omega)$ denotes the weighted Lebesgue space endowed with norm
\[
\|f\|_{L^p(B, \omega)} = \left(\int_B |f(x)|^p d\omega(x) \right)^{\frac{1}{p}}, \quad \text{where} \quad d\omega(x)=\omega(x)dx.
\]
In addition, for $k\in \bN$ and $p \in [1, \infty)$, the Sobolev space $W^{k,p}(\Omega, \omega)$ is defined to be the set of all measurable functions $f$ on $\Omega$ such that their weak derivatives $D^{\beta}f$ are in $L^{p}(\Omega, \omega)$ for all $\beta = (\beta_1,\beta_2,...,\beta_n)\in (\bN \cup \{0\})^n$ satisfying $|\beta|= \beta_1+\beta_2+...+\beta_n\leq k$. The space $W^{k,p}(\Omega, \omega)$ is equipped with the norm
\begin{equation*}
	\|f\|_{W^{k,p}(\Omega,\omega)} =  \sum_{|\beta|\leq k}\|D^{\beta}f\|_{ L^{p}(\Omega, \omega)}.
\end{equation*}
Also, for a measurable set $E\subset\bR^n$, we write 
\[|E|= \int_E dx \quad \text{and} \quad \omega(E)=\int_E \omega(x)dx.\]
For each $\rho>0$ and $x_0 \in \mathbb{R}^n$, we write $B_\rho(x_0)$ for the ball in $\mathbb{R}^n$ centered at $x_0$ with radius $\rho$. For brevity, we write $B_\rho = B_\rho(0)$. Also, for each $x\in\overline{\Omega}$ and $\rho>0$, we denote
\[\Omega_\rho(x)=\Omega\cap B_\rho(x).\]
\smallskip
Next, we provide a few definitions needed in the statement of the main result. The first definition is about the weighted mean oscillation of $\bA$.
\begin{definition}\label{A-bmo-def} Let $\bA:\Omega\times\overline{\bK}\times\bR^n$ satisfy \eqref{bA-growth}, and let $\mu \in L^1_{\text{loc}}(\Omega)$ be nonnegative. For every $x_0\in\overline{\Omega}$ and $\rho>0$, the weighted mean oscillation of $\bA$ in $\Omega_\rho(x_0) = \Omega \cap B_\rho(x_0)$ with respect to the weight $\mu$ is defined by
	\[\Theta_{\rho,x_0}(\bA,\mu) = \dfrac{1}{\mu(\Omega_\rho(x_0))}\int_{\Omega_\rho(x_0)} \Big(\sup_{z\in \bK}\sup_{\xi\in\bR^n\setminus\{0\}}\dfrac{|\bA(x,z,\xi)-\overline{\bA}_{\Omega_\rho(x_0)}(z,\xi)|}{|\xi|}\Big )d\mu(x),\]
	where $\mu(\Omega_\rho(x_0))= \int_{\Omega_\rho(x_0)}\mu(x)dx$, and $\overline{\bA}_{\Omega_\rho(x_0)}(z,\xi)$ denotes the weighted average given by
	\[
	{\overline{\bA}}_{\Omega_\rho(x_0)}(z,\xi) = \dfrac{1}{\mu(\Omega_\rho(x_0))}\int_{\Omega_\rho(x_0)}\bA(x,z,\xi)d\mu(x).\]
\end{definition}
Next, we give the definition on local BMO condition for the class of solutions.
\begin{definition} \label{BMO}
	Let $\rho_0>0$, $f\in L^1_{\text{loc}}(\Omega)$, and let $\omega: \bR^n \rightarrow [0, \infty)$ be locally integrable. We say that $f \in \textup{BMO}_{\rho_0}(\Omega, \omega)$, the class of bounded mean oscillation functions with respect to the weight {$\omega$}, if 
	\[[[f]]_{\textup{BMO}_{\rho_0}(\Omega,\omega)}= \sup_{y\in\Omega}\sup_{0<\rho<\rho_0}\dfrac{1}{\omega(\Omega_{\rho}(y))}\int_{\Omega_{\rho}(y)}|f(x)-\overline{f}_{\Omega_{\rho}(y)}|d\omega(x) <\infty\]
where
\[
\bar{f}_{E} =   \frac{1}{\omega(E)}\int_E f(x)\,d\omega(x),\qquad 0<\omega(E)<\infty.
\] 
\end{definition}
\noindent
{We point out that when $\omega \in A_\infty$,} the equivalence of weighted and unweighted BMO is proved in \cite[Theorem 5]{Muckenhoupt-Weeden-BMO}. From this and with $\omega \in A_\infty$, it follows that a local John-Nirenberg lemma also holds for $f \in \textup{BMO}_{\rho_0}(\Omega, \omega)$, see  \cite[Theorem 4.1]{Buckley-JN} for example. 

\smallskip
Our next definition is about the smoothness of the boundary $\partial \Omega$. Here, $B_\rho'(x')$ be the ball in $\bR^{n-1}$ of radius $\rho>0$ centered at $x' \in \mathbb{R}^{n-1}$.
\begin{definition}\label{lipschitz-domain} 	For given positive numbers $\delta$ and $\rho_0$, we say that $\Omega$ is of $(\delta,\rho_0)$-Lipschitz if for any $x_0=(x_0',x_{n0})\in \p \Omega$ with $x_0'\in \bR^{n-1}$ and $x_{n0}\in\bR$, there exists a Lipschitz continuous function $\gamma:\bR^{n-1}\rightarrow\bR$ such that upon relabelling and reorienting the coordinates, we obtain
	\begin{align*}
		\{x=(x',x_n)\in B_{\rho_0}'(x_0')\times\bR: \gamma(x')<x_n<\gamma(x')+\rho_0\}\subset \Omega,\\
		\{(x',\gamma(x')): x'\in  B_{\rho_0}'(x_0')  \}\subset \p \Omega,
	\end{align*}
	and 
	\[\gamma(x_0')=x_{n0},\quad \nabla \gamma(x_0') = 0,\quad \|\nabla\gamma\|_{L^{\infty}(\bR^{n-1})}\leq \delta.\]
\end{definition} \noindent
{Observe that if $\partial \Omega \in C^1$, then} for every $\delta>0$, there exists some sufficiently small $\rho_0>0$ such that $\Omega$ satisfies the $(\delta,\rho_0)$-Lipschitz condition.	

\smallskip
Now, we give the concept of weak solutions to the {conormal} problem \eqref{main-eqs}.
\begin{definition}\label{weak-sol-def} Let $\bK \subset \bR$ be an interval, $\alpha\in(-1,\infty)$, and $\Omega \subset \bR^n$ be open and bounded with Lipschitz boundary $\p \Omega$. Assume that $\bA: \Omega\times\overline{\bK}\times\bR^n \rightarrow \bR^n$ satisfies \eqref{bA-growth}. Let $\bF\in L^p(\Omega,\mu)^n$ for some $1<p<\infty$, and assume that $\mu$ satisfies \eqref{mu-weight}. A function $u\in W^{1,p}(\Omega,\mu)$ is said to be a weak solution of \eqref{main-eqs} {if $u(x)\in\overline{\bK}$ for almost every $x\in\Omega$ and }
	\begin{equation}
		\int_{\Omega}\mu(x)\langle \bA(x, u(x),\nabla u(x)), \nabla\varphi(x)\rangle dx = \int_{\Omega}\mu(x)\langle \bF(x),\nabla \varphi(x) \rangle  dx,
	\end{equation}
	for every $\varphi\in C^{\infty}(\overline{\Omega})$.
\end{definition}
\smallskip
\noindent
In the statement of Theorem \ref{main-thrm} below, for a given $r_0$ in \eqref{mu-weight} and a small  constant $\rho_0 \in (0,1)$  representing the Lipschitz geometry of $\Omega$, we denote $\bar{\rho} = \min\{r_0, \rho_0\}/32$ and
\[
\Omega^{\bar{\rho}} = \{x \in \Omega:   \operatorname{dist}(x,\partial\Omega) > \rho \}.
\]
In addition,  the two positive constants $c_0, C_0$ depending on $\mu, r_0$ and $\rho_0$ are defined by
\begin{equation} \label{bound-mu}
0<c_0\le\mu(x)\le C_0<\infty
\quad\text{if}\quad
\bar{\rho}/2<\operatorname{dist}(x,\partial\Omega)\le2r_0.
\end{equation}
We now state the main result. 
\begin{theorem} \label{main-thrm}  Let $\Lambda\geq1$, $M>0$, $\alpha\in(-1,\infty)$,  and $p\in [2,\infty)$. There is a sufficiently small positive number $\delta = \delta(\Lambda, \alpha, p, M, n)$ such that the following holds. Suppose that $\bK \subset \bR$ is an open interval, $\partial \Omega$ is $(\delta, \rho_0)$-Lipschitz with some $\rho_0 \in (0,1)$. Suppose also that $\bA:\Omega\times\overline{\bK}\times {\bR^n} \rightarrow \bR^n$ is a vector field satisfying \eqref{bA-growth}, $\mu$ satisfies \eqref{mu-weight}, and 
	\begin{equation}\label{A-bmo-main}
		\sup_{\rho\in(0,R_0)}\sup_{x\in \overline{\Omega}}\Theta_{\rho,x}(\bA,\mu) \leq \delta,
\end{equation}
for some $R_0\in(0,1)$. Then for any weak solution $u\in W^{1,2}(\Omega,\mu)$ of \eqref{main-eqs} satisfying $[[u]]_{\textup{BMO}_{R_0}(\Omega,\mu)}\leq M$ and with $\bF\in L^p(\Omega,\mu)^n$,  we have $u\in W^{1,p}(\Omega,\mu)$ and 
\begin{equation} \label{main-estimate2}
\|u - \bar{u}_{\Omega}\|_{{W^{1,p}}(\Omega,\mu)} \leq  C\|\bF\|_{L^{p}(\Omega,\mu)},
\end{equation}
where  $C=C(\Lambda,M,\alpha,p,n,\bK, {R_0,\rho_0,r_0,\Omega,c_0,C_0})>0$ also depends on the modulus of continuity of $\mu$ on $\overline{\Omega^{\bar{\rho}}}$ with $\bar{\rho} = \min\{\rho_0, r_0\}/32$ and
\[
\bar{u}_{\Omega} = \frac{1}{\mu(\Omega)} \int_{\Omega} u(x) d\mu(x).
\]
\end{theorem}

\smallskip
The proof of Theorem~\ref{main-thrm} is based on a perturbation argument and weighted level-set estimates. We first establish weighted Poincar\'{e} and Sobolev-Poincar\'{e} inequalities for powers of the distance to the boundary with $\alpha>-1$, including weights outside the $A_2$ class. These inequalities are used to derive H\"older estimates and higher integrability for solutions of the homogeneous comparison equations. We then construct comparison solutions by replacing the weight by its flat model, freezing the dependence of the coefficients on the solution, and averaging the coefficients in the spatial variable. The weighted BMO bound on $u$ and the small mean oscillations of the coefficients control the approximation errors. Combining the comparison estimates with weighted maximal-function level-set estimates yields local $L^p$-estimates for the gradient in flat boundary cylinders. Finally, we flatten the boundary and combine the boundary and interior estimates with the energy inequality and the weighted Sobolev-Poincar\'{e} inequality to obtain \eqref{main-estimate2}. To deal with the general setting in which $\bA$ depends on both $u$ and $\nabla u$,  we adapt a scaling constant technique developed in \cite{NP, P}  to carry out the normalizations in the proof. As such, we introduce the parameter $\lambda$ and consider the enlarged one-parameter class of equations:
\begin{equation}\label{main-eqs-lambda}
	\left\{
	\begin{array}{cccl}
		\diver(\mu(x)\bA(x,\lambda u,\nabla u)) &= &\diver(\mu(x)\bF(x)) &\text{in }\Omega,\\
		\mu(x)\Big(\bA(x,\lambda u,\nabla u)-\bF(x)\Big)\cdot\overrightarrow\nu(x) &= &0 &x \in\partial\Omega,
	\end{array}\right.
\end{equation}
with $\lambda\geq0$.  The important structural feature in \eqref{main-eqs-lambda} that enables us to apply the perturbation technique in \cite{CP} is that this class of equations is invariant under the scaling $u(x) \mapsto u(r x)/r$ for all $r>0$. The original equation \eqref{main-eqs} corresponds to $\lambda=1$ and it is a subclass of \eqref{main-eqs-lambda}. The case that $\bA$ is independent on $u$ as in \cite{Dong-Phan-Sire} can be considered as $\lambda=0$ and $\bK$ is any interval containing $0$. 

\smallskip
The paper is organized as follows. In Section~\ref{sec:FA-Results}, we recall the necessary definitions and develop preliminary results. 
 In Section~\ref{sec:Lp-flat}, we establish local boundary $L^p$-estimates in flat domains. Finally, in Section~\ref{sec:main-proofs}, we derive the interior and boundary gradient estimates and prove Theorem~\ref{main-thrm}. 
\section{Muckenhoupt weights, Weighted Poincar\'{e} inequalities, and preliminary estimates}\label{sec:FA-Results}
\subsection{Muckenhoupt weights} 
We begin by recalling the definition of the $A_p$-Muckenhoupt class of weights.
\begin{definition}
	Let $p\in (1,\infty)$. A nonnegative locally integrable function $\omega:\bR^n\rightarrow [0,\infty)$ is called an $A_p$ weight, if the quantity
	\[[\omega]_p =\sup_{\text{balls }B\subset \bR^n} \left(\dfrac{1}{|B|}\int_{B}\omega(x)dx\right)\left(\dfrac{1}{|B|}\int_{B}\omega(x)^{\frac{1}{1-p}}dx\right)^{p-1} < \infty.\]
\end{definition}
\noindent It follows from the definition that if $\omega \in A_p$, then $[\omega]_{p} \geq 1$. We use the following properties of $A_p$ weights; see \cite[Proposition 9.1.5 and Theorem 9.2.2]{Grafakos} and \cite[Theorem IV]{CF}.
\begin{lemma}\label{lm-doubling} For $1<p<\infty$, the following statements hold 
\begin{itemize}
\item[\textup{(i)}] For every $\omega \in A_p$, 
\[
\frac{\omega(B)}{\omega(E)}\leq [\omega]_p \Big(\dfrac{|B|}{|E|}\Big)^p
\]
 for every ball $B \subset \bR^n$, and for every measurable set $E \subset B$ with $|E| \not=0$.
\item[\textup{(ii)}] For $M_0 \geq 1$, there are $C=C(M_0,p, n)>0$ and $\vartheta=\vartheta(M_0,p, n) \in (0,1)$ so that
\[C \Big(\dfrac{|B|}{|E|}\Big)^\vartheta\leq  \frac{\omega(B)}{\omega(E)}\]
for every $\omega \in A_p$ satisfying $[\omega]_p\leq M_0$, and for every measurable set $E$  and for every ball $B\subset\bR^n$ such that $E \subset B$ and $|E| \not=0$.
\end{itemize}
\end{lemma}

\begin{remark}\label{weight-remark} Let $\alpha \in (-1, \infty)$ and let $\omega: \mathbb{R}^{n-1}\times \bR \rightarrow \bR$ be defined by 
	\begin{equation}
		\omega(x) = |x_{n}|^{\alpha},  \quad \forall \  x = (x',x_n)\in \mathbb{R}^{n-1} \times \mathbb{R}.
	\end{equation}
Then, $\omega\in A_p$ for every $p>\max\{1,\alpha+1\}$, with $[\omega]_p\le C(n,\alpha,p)$. In particular,
\[
\omega\in A_{2+\max\{\alpha,0\}},\qquad
[\omega]_{2+\max\{\alpha,0\}}\le M_0(n,\alpha).
\]
\end{remark}



\subsection{Weighted Poincar\'{e} and Sobolev-Poincar\'{e} inequalities} \label{Poincare-ineq}
Weighted Poincar\'{e} inequalities for $A_2$ weights are well known; see \cite{FKS}. In the setting of this paper, we need the following weighted Poincar\'{e} inequality with the weight $x_n^\alpha$, $\alpha>-1$, which need not belong to $A_2$ when $\alpha\ge1$.
\begin{lemma} \label{weighted-poincare-top}
Let $R > 0$, let $D^+_R = B'_R \times (0, R) \subset \mathbb{R}^n$ be the half-cylinder where $B'_R = \{x' \in \mathbb{R}^{n-1} : |x'| < R\}$, and let $\alpha > -1$. Then the weighted Poincar\'{e} inequality
\[
\int_{D^+_R} |u(x)|^2 x_n^\alpha \, dx \le \frac{4 R^2}{(\alpha+1)^2}  \int_{D^+_R} |\partial_n u(x)|^2 x_n^\alpha \, dx
\]
holds for every $u \in C^1(\overline{D_R^+})$ that vanishes on $B'_R \times \{R\}$. 
\end{lemma}
\begin{proof} Let us fix $x' \in B'_R$ and define $f(x_n) = u(x', x_n)$ for $x_n \in (0, R)$. By the assumption that $u(x', R)= 0$, we have $f(R) = 0$. Then, it follows from the fundamental theorem of calculus that 
\[
f(x_n) = -\int_{x_n}^R f'(t) \, dt \quad \text{for any} \quad x_n \in (0, R).
\]
Let us denote  the auxiliary function $F(x_n) = \int_{x_n}^R f'(t) \, dt = -f(x_n)$. Then, we see from an integration by parts and $\alpha +1>0$ that
\begin{align*}
\int_0^R |F(x_n)|^2 x_n^\alpha \, dx_n & = - \frac{2}{\alpha+1}\int_0^R F(x_n) F'(x_n)x_n^{\alpha+1} \, dx_n.
\end{align*}
Substituting $F'(x_n) = -f'(x_n)$ into this last equality, we obtain
\[
\int_0^R |f(x_n)|^2 x_n^\alpha \, dx_n = \frac{2}{\alpha+1} \int_0^R F(x_n) f'(x_n) x_n^{\alpha+1} \, dx_n.
\]
Note that by applying Cauchy-Schwarz to the right-hand side, we infer that
\[
\left| \int_0^R F(x_n) f'(x_n) x_n^{\alpha+1} \, dx_n \right| \le \left( \int_0^R |F(x_n)|^2 x_n^\alpha \, dx_n \right)^{1/2} \left( \int_0^R |f'(x_n)|^2 x_n^{\alpha+2} \, dx_n \right)^{1/2}.
\]
Therefore, it follows that
\[
\int_0^R |f(x_n)|^2 x_n^\alpha \, dx_n \le \frac{4}{(\alpha+1)^2} \int_0^R |f'(x_n)|^2 x_n^{\alpha+2} \, dx_n.
\]
Since $x_n\le R$, we have
\[
x_n^{\alpha+2} = x_n^2 x_n^\alpha \le R^2 x_n^\alpha.
\]
Hence,
\[
\int_0^R |u(x', x_n)|^2 x_n^\alpha \, dx_n \le \frac{4 R^2}{(\alpha+1)^2} \int_0^R |\partial_n u(x', x_n)|^2 x_n^\alpha \, dx_n.
\]
Finally, integrating over $B'_R$, we obtain
\[
\int_{D^+_R} |u(x)|^2 x_n^\alpha \, dx \le \frac{4 R^2}{(\alpha+1)^2} \int_{D^+_R} |\partial_n u(x)|^2 x_n^\alpha \, dx.
\]
This completes the proof.
\end{proof}

\smallskip
In various places, we need the following weighted Poincar\'e inequality for functions that do not vanish on parts of the boundary.
\begin{lemma} \label{weighted-poincare-mean}
Let $R>0$, let $D_R^+=B'_R\times(0,R)$ where
$B'_R=\{x'\in\mathbb R^{n-1}:|x'|<R\}$, and let
$\alpha>-1$. Then, there is $C=C(n,\alpha)>0$
such that the weighted Poincar\'{e} inequality
\[
\int_{D_R^+}|u(x)-\bar u|^2x_n^\alpha\,dx
\le CR^2\int_{D_R^+}|\nabla u(x)|^2x_n^\alpha\,dx
\]
holds for any $u\in C^1(\overline{D_R^+})$,
where
\[
\bar u=
\left(\int_{D_R^+}x_n^\alpha\,dx\right)^{-1}
\int_{D_R^+}u(x)x_n^\alpha\,dx.
\]
\end{lemma}
\begin{proof} For any fixed $x'\in B'_R$, let
$f(t)=u(x',t)$ for $t\in(0,R)$, 
\[
m=\int_0^R t^\alpha\,dt
=\frac{R^{\alpha+1}}{\alpha+1},
\quad\text{and}\quad
f_w(x')=\frac1m\int_0^R f(t)t^\alpha\,dt.
\]
Since $f_w$ is the weighted mean of $f$, for any
constant $c$, we have
\[
\int_0^R|f(t)-c|^2t^\alpha\,dt
=\int_0^R|f(t)-f_w|^2t^\alpha\,dt +m|f_w-c|^2,
\]
where we used
\[
\int_0^R(f(t)-f_w)t^\alpha\,dt=0.
\]
Taking $c=f(R)$, we obtain
\begin{align} \notag
\int_0^R|f(t)-f_w|^2t^\alpha\,dt  & = \int_0^R|f(t)-f(R)|^2t^\alpha\,dt - m|f_w-f(R)|^2\\ \label{wei-mean-ineq}
& \le\int_0^R|f(t)-f(R)|^2t^\alpha\,dt.
\end{align}
To estimate the right-hand side of this last inequality, let us denote
$h(t)=f(t)-f(R)$. Since $h(R)=0$ and $\alpha>-1$, we apply the
integration by parts and the Cauchy-Schwarz
inequality to infer
\[
\begin{aligned}
\int_0^R|h(t)|^2t^\alpha\,dt
&=-\frac{2}{\alpha+1}
\int_0^R h(t)h'(t)t^{\alpha+1}\,dt \le\frac{2}{\alpha+1}
\left(\int_0^R|h|^2t^\alpha\,dt\right)^{1/2}
\left(\int_0^R|h'|^2t^{\alpha+2}\,dt\right)^{1/2}\\
&\le\frac{2R}{\alpha+1}
\left(\int_0^R|h|^2t^\alpha\,dt\right)^{1/2}
\left(\int_0^R|h'|^2t^\alpha\,dt\right)^{1/2}
\end{aligned}
\]
Therefore,
\[
\left(\int_0^R|h(t)|^2t^\alpha\,dt \right)^{1/2} \leq  \frac{2R}{\alpha+1} \left(\int_0^R|h'|^2t^\alpha\,dt\right)^{1/2}
\]
and consequently,
\[
\int_0^R|f(t)-f_w|^2t^\alpha\,dt
\le\frac{4R^2}{(\alpha+1)^2}
\int_0^R|f'(t)|^2t^\alpha\,dt.
\]
Integrating this inequality over $x'\in B'_R$, we obtain
\begin{equation} \label{1D-ineq-Nov-4}
\int_{D_R^+}|u(x',x_n)-f_w(x')|^2x_n^\alpha\,dx
\le\frac{4R^2}{(\alpha+1)^2}
\int_{D_R^+}|\partial_nu|^2x_n^\alpha\,dx.
\end{equation}

\smallskip
Now, note that
\[
\bar u=\frac1{|B'_R|}
\int_{B'_R}f_w(x')\,dx'
\]
and 
\[
\begin{aligned}
m|\nabla'f_w(x')|^2
&=\frac1m\left|\int_0^R\nabla'u(x',t)t^\alpha\,dt\right|^2 \le\frac1m
\left(\int_0^R|\nabla'u(x',t)|^2t^\alpha\,dt\right)
\left(\int_0^R t^\alpha\,dt\right)\\
&=\int_0^R|\nabla'u(x',t)|^2t^\alpha\,dt.
\end{aligned}
\]
Then, applying the Poincar\'{e} inequality on
$B'_R$, we obtain
\begin{equation} \label{Poincare-x'-Nov-4}
m\int_{B'_R}|f_w-\bar u|^2\,dx' \le C_nR^2m
\int_{B'_R}|\nabla'f_w|^2\,dx' \le C_nR^2
\int_{D_R^+}|\nabla'u(x)|^2x_n^\alpha\,dx.
\end{equation}
Finally, using the fact that
\[
\int_0^R(u(x',t)-f_w(x'))t^\alpha\,dt=0,
\]
we infer from \eqref{1D-ineq-Nov-4} and \eqref{Poincare-x'-Nov-4} that
\[
\begin{aligned}
\int_{D_R^+}|u-\bar u|^2x_n^\alpha\,dx
&=\int_{D_R^+}|u-f_w|^2x_n^\alpha\,dx
+m\int_{B'_R}|f_w-\bar u|^2\,dx'\\
&\le CR^2
\int_{D_R^+}|\nabla u(x)|^2x_n^\alpha\,dx,
\end{aligned}
\]
where $C=C(n,\alpha)>0$. This completes the proof.
\end{proof}

\smallskip
The following weighted Sobolev-Poincar\'e inequality follows from \cite[Theorem 9.7, pp.~43-44]{HK-Poincare} and \cite[Remark 3.2]{Dong-Phan}. Here, for completeness, we provide a direct proof.
\begin{lemma} \label{weighted-sobolev-poincare}
Let $R>0$, and $D_R^+=B'_R\times(0,R)$ where
$B'_R=\{x'\in\mathbb R^{n-1}:|x'|<R\}$.
Also, let $\alpha>-1$, $N=n+\alpha$,
$Q=n+\max\{\alpha,0\}$, and $q$ satisfy
\[
q=\begin{cases}
\frac{2Q}{Q-2} & \text{if } Q>2,\\
\text{any }q\in[2,\infty) & \text{if } Q\le2.
\end{cases}
\]
Then, there is $C=C(n,q,\alpha)>0$ such that
the weighted Sobolev-Poincar\'{e} inequality
\[
\left(\int_{D_R^+}|u(x)-\bar u|^q x_n^\alpha\,dx
\right)^{2/q}
\le CR^{2-N(1-2/q)}
\int_{D_R^+}|\nabla u(x)|^2x_n^\alpha\,dx
\]
holds for any $u\in C^1(\overline{D_R^+})$,
where
\[
\bar u=
\left(\int_{D_R^+}x_n^\alpha\,dx\right)^{-1}
\int_{D_R^+}u(x)x_n^\alpha\,dx.
\]
\end{lemma}
\begin{proof}
By scaling, we may assume without loss of generality that $R=1$. For $q=2$, we apply the preceding lemma. Hence, we only need to consider the case $q>2$. Let us denote $s=2q/(q-2)$ and note that $s\ge Q$ and $s>2$. Let us also denote $d\mu(y)=y_n^\alpha\,dy$. For
$x\in D_1^+$ and $0<r\le2$, let
\begin{equation} \label{E_r-u_r}
\begin{split}
& E_r(x)=D_1^+\cap
\bigl(B'_r(x')\times(x_n-r,x_n+r)\bigr),
\quad \text{and}\\
& u_r(x)=\frac{1}{\mu(E_r(x))}
\int_{E_r(x)}u(y)\,d\mu(y).
\end{split}
\end{equation}
Note that $E_2(x)=D_1^+$ and $u_2(x)=\bar u$. Also, it follows from a direct calculation
that there are constants $c = c(n, \alpha)>0$ and $C= C(n, \alpha)>0$ such that
\[
cr^n(x_n+r)^\alpha
\le\mu(E_r(x))
\le Cr^n(x_n+r)^\alpha,
\qquad 0<r\le2.
\]
Note that if $\alpha\ge0$, then
$(x_n+r)^\alpha\ge r^\alpha$, while if
$-1<\alpha<0$, we have
$(x_n+r)^\alpha\ge3^\alpha$.
Consequently,
\begin{equation} \label{lower-mu-E_r}
\mu(E_r(x))\ge {Cr^Q\ge C\,2^{Q-s}r^s},
\qquad 0<r\le2.
\end{equation}
Furthermore,
\begin{equation} \label{mu-doubling}
\mu(E_{2r}(x))\le C\mu(E_r(x)),
\qquad 0<r\le1.
\end{equation}
Here, the constants $C = C(n, \alpha, q)>0$ and can change line by line.

\smallskip
We next establish a local weighted Poincar\'{e} inequality. For $I=(a,b)\subset(0,1)$, with
$a=0$ also allowed, let
\[
f_{I,\alpha}=\frac1{m_I}
\int_I f(t)t^\alpha\,dt \quad \text{with} \quad m_I=\int_I t^\alpha\,dt.
\]
Since $\alpha>-1$, we have
\[
\int_a^t\tau^\alpha\,d\tau
\le C_\alpha(b-a)t^\alpha,
\qquad a<t<b.
\]
Using the fact that $f_{I,\alpha}$ is the
weighted mean of $f$ as in \eqref{wei-mean-ineq}, we obtain
\begin{align} \notag
\int_I|f-f_{I,\alpha}|t^\alpha\,dt
&\le2\int_I|f(t)-f(b)|t^\alpha\,dt \le2\int_a^b|f'(t)|
\left(\int_a^t\tau^\alpha\,d\tau\right)dt
\notag\\
&\le C_\alpha(b-a)
\int_I|f'(t)|t^\alpha\,dt.
\label{enq5-Nov-4}
\end{align}
The estimate also holds when $a=0$ by
applying the fundamental theorem of calculus
on intervals away from zero and passing
to the limit.

\smallskip
Now, for fixed $x\in D_1^+$ and $0<r\le2$,
let
\[
\Omega=B'_1\cap B'_r(x'),
\quad
I=(0,1)\cap(x_n-r,x_n+r),
\quad\text{and}\quad
v(y')=\frac1{m_I}
\int_Iu(y',t)t^\alpha\,dt.
\]
Note from \eqref{E_r-u_r} that
\[
u_r(x)=\frac1{|\Omega|}
\int_\Omega v(y')\,dy'.
\]
Applying \eqref{enq5-Nov-4} to $u(y',\cdot)$ on $I$
and integrating over $y'\in\Omega$, we obtain
\begin{equation}\label{ponce-4-Nove-4}
\int_\Omega\int_I
|u(y',t)-v(y')|t^\alpha\,dt\,dy'
\le Cr\int_{E_r(x)}
|\partial_nu(y)|\,d\mu(y).
\end{equation}
Since $\Omega$ is convex and has diameter
at most $2r$, it follows from the Poincar\'{e}
inequality that
\begin{equation}\label{ponce-5-Nove-4}
m_I\int_\Omega|v(y')-u_r(x)|\,dy'
\le Crm_I\int_\Omega|\nabla'v(y')|\,dy'\\
\le Cr\int_{E_r(x)}
|\nabla'u(y)|\,d\mu(y),
\end{equation}
where we used
\[
m_I|\nabla'v(y')|
\le\int_I|\nabla'u(y',t)|t^\alpha\,dt.
\]
Combining \eqref{ponce-4-Nove-4} and
\eqref{ponce-5-Nove-4}, we obtain
\begin{align}\notag
\int_{E_r(x)}|u(y)-u_r(x)|\,d\mu(y)
&\le\int_\Omega\int_I
|u(y',t)-v(y')|t^\alpha\,dt\,dy'
+m_I\int_\Omega|v(y')-u_r(x)|\,dy'\\ \notag
&\le Cr\int_{E_r(x)}
\bigl(|\partial_nu(y)|+|\nabla'u(y)|\bigr)
\,d\mu(y)\\ \label{local-poincare-Nov-4}
&\le Cr\int_{E_r(x)}
|\nabla u(y)|\,d\mu(y).
\end{align}

\smallskip
Now, let us recall the Hardy-Littlewood maximal operator associated with the weighted measure $d\mu(y)=y_n^\alpha\,dy$ and the sets $E_r(x)$:
\[
\mathcal M_\mu h(x)
=\sup_{0<r\le2}
\frac{1}{\mu(E_r(x))}
\int_{E_r(x)}|h(y)|\,d\mu(y)
\]
for any $h\in L^1(D_1^+,d\mu)$. Note that, by the doubling property \eqref{mu-doubling}, there exists $C=C(n,\alpha)>0$ such that the following weak-type $(1,1)$ estimate holds:
\begin{equation} \label{weak-type-est-HL}
\mu\bigl(\{x\in D_1^+:\mathcal M_\mu h(x)>t\}\bigr)
\le\frac Ct\int_{D_1^+}|h(y)|\,d\mu(y)
\end{equation}
holds for every \(t>0\).

\smallskip
Next, let us denote
\begin{equation} \label{G-A-def-nov-4}
\begin{split}
& G(x)=\sup_{0<r\le2}
\frac{1}{\mu(E_r(x))}\int_{E_r(x)}|\nabla u(y)|\,d\mu(y),
\quad \text{and} \\
& A=\left(\int_{D_1^+}|\nabla u(x)|^2
\,d\mu(x)\right)^{1/2}.
\end{split}
\end{equation}
We first estimate $G$ in $L^2(D_1^+,d\mu)$.  For any $t>0$, by writing
\[
|\nabla u|=g_1+g_2,
\quad \text{where} \quad
g_1=|\nabla u|\mathbf 1_{\{|\nabla u|>t/2\}},
\qquad
g_2=|\nabla u|\mathbf 1_{\{|\nabla u|\le t/2\}},
\]
and applying the weak-type estimate \eqref{weak-type-est-HL} to $g_1$,
we obtain
\[
\mu(\{G>t\})
\le\frac Ct
\int_{\{|\nabla u|>t/2\}}
|\nabla u(x)|\,d\mu(x).
\]
Then, integrating in $t$ and applying Fubini's
theorem, we infer
\begin{align} \notag
\int_{D_1^+}G(x)^2\,d\mu(x)
&=2\int_0^\infty t\,\mu(\{G>t\})\,dt \le C\int_0^\infty
\int_{\{|\nabla u|>t/2\}}
|\nabla u(x)|\,d\mu(x)\,dt\\ \notag
&=C\int_{D_1^+}|\nabla u(x)|
\left(\int_0^{2|\nabla u(x)|}dt\right)d\mu(x)\\ \label{G-A-compare}
&\le C\int_{D_1^+}|\nabla u(x)|^2\,d\mu(x)
=CA^2.
\end{align}

\smallskip
Next, let $r_j=2^{1-j}$, $j\ge0$. From \eqref{E_r-u_r} and the Lebesgue differentiation theorem,
$u_{r_j}(x)\to u(x)$ for almost every
$x\in D_1^+$. Since $u_{r_0}(x)=\bar u$,
we have
\[
|u(x)-\bar u|
\le\sum_{j=0}^\infty
|u_{r_{j+1}}(x)-u_{r_j}(x)|.
\]
Using \eqref{local-poincare-Nov-4} and
the doubling estimate \eqref{mu-doubling}, we obtain
\[
\begin{aligned}
|u_{r_{j+1}}(x)-u_{r_j}(x)|
&\le
\frac1{\mu(E_{r_{j+1}}(x))}
\int_{E_{r_{j+1}}(x)}
|u(y)-u_{r_j}(x)|\,d\mu(y)\\
&\le
\frac C{\mu(E_{r_j}(x))}
\int_{E_{r_j}(x)}
|u(y)-u_{r_j}(x)|\,d\mu(y)\\
&\le  \frac{Cr_j}{\mu(E_{r_j}(x))} \int_{E_{r_j}(x)}
|\nabla u(y)|\,d\mu(y).
\end{aligned}
\]
Consequently,
\begin{equation} \label{u-baru-est-Nov-4}
|u(x)-\bar u|
\le C\sum_{j=0}^\infty 
\frac{r_j}{\mu(E_{r_j}(x))} \int_{E_{r_j}(x)}
|\nabla u(y)|\,d\mu(y).
\end{equation}
We now estimate the sum on the right-hand side of \eqref{u-baru-est-Nov-4}. Let $0<\rho\le1$ be determined later, and split the sum into
$r_j\le\rho$ and $r_j>\rho$.
By the definition of $G(x)$, we note that
\begin{equation} \label{smaller-term-Nov-4}
\begin{aligned}
\sum_{r_j\le\rho}
\frac{r_j}{\mu(E_{r_j}(x))} \int_{E_{r_j}(x)}
|\nabla u(y)|\,d\mu(y)
&\le G(x)\sum_{r_j\le\rho}r_j \le C\rho G(x).
\end{aligned}
\end{equation}
On the other hand, by H\"older's inequality
and the lower bound for $\mu(E_{r_j}(x))$ in \eqref{lower-mu-E_r}, it follows that
\[
\begin{aligned}
\frac{1}{\mu(E_{r_j}(x))}\int_{E_{r_j}(x)}
|\nabla u(y)|\,d\mu(y)
&\le\frac1{\mu(E_{r_j}(x))^{1/2}}
\left(\int_{E_{r_j}(x)}
|\nabla u(y)|^2\,d\mu(y)\right)^{1/2}\\
&\le CA r_j^{-s/2}.
\end{aligned}
\]
From this and as $s>2$, we infer that
\[
\begin{aligned}
\sum_{r_j>\rho}
\frac{r_j}{\mu(E_{r_j}(x))} \int_{E_{r_j}(x)}
|\nabla u(y)|\,d\mu(y)
&\le CA\sum_{r_j>\rho}r_j^{1-s/2}\le CA\rho^{1-s/2}.
\end{aligned}
\]
Then, by plugging \eqref{smaller-term-Nov-4} and the last estimate into \eqref{u-baru-est-Nov-4}, we obtain
\[
|u(x)-\bar u|
\le C\rho G(x)+CA\rho^{1-s/2},
\qquad 0<\rho\le1.
\]
If $A=0$, then it follows from \eqref{G-A-def-nov-4} that $u=\bar u$ almost everywhere and the assertion of the lemma is trivial.
Otherwise, $G(x)>0$, and we take
\[
\rho=\min\left\{1,
\left(\frac A{G(x)}\right)^{2/s}\right\}.
\]
If $G(x)\ge A$, then
\[
\rho G(x)=A^{2/s}G(x)^{1-2/s},
\quad \text{and} \quad
A\rho^{1-s/2}=A^{2/s}G(x)^{1-2/s}.
\]
If $G(x)<A$, then $\rho=1$ and
$\rho G(x)+A\rho^{1-s/2}\le2A$.
Thus, in either case, we obtain
\[
|u(x)-\bar u|
\le C\left(A+A^{2/s}G(x)^{1-2/s}\right).
\]
Since $q(1-2/s)=2$ and $2q/s=q-2$,
we have
\[
\begin{aligned}
\int_{D_1^+}|u(x)-\bar u|^q\,d\mu(x)
&\le CA^q+
CA^{2q/s}\int_{D_1^+}
G(x)^{q(1-2/s)}\,d\mu(x)\\
&=CA^q+CA^{q-2}
\int_{D_1^+}G(x)^2\,d\mu(x)\\
&\le CA^q,
\end{aligned}
\]
where we use \eqref{G-A-compare} in the last estimate. This proves the assertion of the lemma when $R=1$. By making the change of variables
$x=Ry$, we obtain
\[
\left(\int_{D_R^+}|u(x)-\bar u|^q
x_n^\alpha\,dx\right)^{2/q}
\le CR^{2-N(1-2/q)}
\int_{D_R^+}|\nabla u(x)|^2x_n^\alpha\,dx.
\]
The proof is complete.
\end{proof}

\begin{remark} When $\alpha\ge1$, the Poincar\'{e} type inequality in Lemma \ref{weighted-poincare-mean} on $B_1$ fails. Indeed, by \cite{STV-even}, the function $u(x)=\operatorname{sgn}(x_n)$  satisfies
\[
\bar u=0,\qquad
\int_{B_1}|u(x)-\bar u|^2|x_n|^\alpha\,dx>0,
\qquad
\int_{B_1}|\nabla u(x)|^2|x_n|^\alpha\,dx=0.
\]
Note that in this example, the gradient of $u$ is understood in the completion of smooth functions in the norm $L^2(B_1, |x_n|^\alpha)$. Our domain $D_R^+$ lies on one side of $\{x_n=0\}$, where this function is constant. Thus, this example does not contradict the preceding inequalities. 
For compactly supported functions, Sire, Terracini and Vita also prove in \cite{STV-even} a weighted Sobolev inequality with the same effective dimension $Q=n+\max\{\alpha,0\}$ in our notation. See also \cite[Lemma 3.1]{Dong-Phan} for similar results on weighted Sobolev inequality in the parabolic case.
\end{remark}

\smallskip
Weighted Sobolev-Poincar\'{e} inequalities in bounded domains involving powers of the distance to the boundary have been studied extensively; see, for example, \cite{Chua1995,ChuaWheeden2010,DrelichmanDuran2008} and the references therein, and \cite[Lemma 3.1]{Dong-Phan} for the parabolic case. For completeness, we establish and prove the following proposition on weighted Sobolev-Poincar\'{e} inequality for the weight defined in \eqref{mu-weight}. The proof uses the weighted Sobolev-Poincar\'{e} inequality on half-cylinders (Lemma \ref{weighted-sobolev-poincare}), a flattening of the boundary, and a finite partition of unity. Also, in the statement of the proposition, the constants $c_0$ and $C_0$ are defined in  \eqref{bound-mu}.
\begin{proposition} \label{weighted-sobolev-poincare-lipschitz}
Let $\Omega\subset\mathbb R^n$, $n \in \mathbb{N}$, be a
bounded connected Lipschitz domain. Let $\alpha>-1$, $N=n+\alpha$,
$Q=n+\max\{\alpha,0\}$, and $q$ satisfy
\[
q=\begin{cases}
\frac{2Q}{Q-2} & \text{if }Q>2,\\
\text{any }q\in[2,\infty) & \text{if }Q\le2.
\end{cases}
\]
Then, there is $C>0$, depending only on
$n,\alpha,q,r_0,c_0,C_0$, and the Lipschitz
geometry of $\Omega$, such that
\[
\left(\int_\Omega|u(x)-\bar u|^q\mu(x)\,dx
\right)^{2/q}
\le CR^{2-N(1-2/q)}
\int_\Omega|\nabla u(x)|^2\mu(x)\,dx
\]
holds for any $u\in C^1(\overline{\Omega})$,
where $R=\operatorname{diam}(\Omega)$ and
\[
\bar u=
\left(\int_\Omega\mu(x)\,dx\right)^{-1}
\int_\Omega u(x)\mu(x)\,dx.
\]
\end{proposition}
\begin{proof} It is sufficient to prove the proposition when $n \geq 2$. We first prove the inequality when
$\operatorname{diam}(\Omega)=1$.
Let us denote
\[
\delta(x)=\operatorname{dist}(x,\partial\Omega),
\quad \text{and} \quad d\mu(x)=\mu(x)\,dx.
\]
Since $\Omega$ is a bounded Lipschitz domain,
we may choose finitely many boundary patches in $\{x \in \Omega: \delta(x) <r_0\}$
whose union contains $\{x \in \Omega: \delta(x) <2\bar{\rho}\}$, where $\bar{\rho}>0$ is defined in \eqref{bound-mu} with $\rho_0 \in (0,1)$ being the constant representing the Lipschitz geometry of $\Omega$.
We cover the compact set $\{x\in\Omega:\delta(x)\ge2\bar{\rho}\}$
by finitely many interior balls contained in $\{ x \in \Omega: \delta(x)>\bar{\rho}\}$.
Denote the resulting cover by $\{U_i\}_{i=1}^m$.

\smallskip
For each boundary patch $U_i$, after a
translation and rotation, we may write
\[
U_i=\Phi_i(D_{r_i}^+),
\quad \text{and} \quad D_{r_i}^+=B'_{r_i}\times(0,r_i),
\]
where
\[
\Phi_i(y',y_n)=(y',y_n+\varphi_i(y'))
\]
and $\varphi_i : \bR^{n-1} \rightarrow \bR$ is a Lipschitz function satisfying $\varphi_i (0) =0$.  We note that, after a translation and rotation, we may assume
that the boundary point is the origin and
$\varphi_i(0)=0$. Let $L$ be the Lipschitz
constant of $\varphi_i$. Since $\Omega$ is a
Lipschitz domain, there is $a>0$ such that
\[
\Omega\cap B_{4a}
=\{x\in B_{4a}:x_n>\varphi_i(x')\}.
\]
We choose $r_i>0$ sufficiently small so that
\[
r_i\sqrt{1+(1+L)^2}\le a.
\]
Then, for any $y\in D_{r_i}^+$, we have
$\Phi_i(y)\in B_a$. Moreover, for any
$z\in\partial\Omega\cap B_{4a}$, writing
$z=(z',\varphi_i(z'))$, we obtain
\[
\begin{aligned}
y_n
&\le |y_n+\varphi_i(y')-\varphi_i(z')|
+|\varphi_i(z')-\varphi_i(y')|\\
&\le |y_n+\varphi_i(y')-\varphi_i(z')|
+L|y'-z'|\\
&\le \sqrt{1+L^2}\,|\Phi_i(y)-z|.
\end{aligned}
\]
On the other hand, if
$z\in\partial\Omega\setminus B_{4a}$, then
\[
|\Phi_i(y)-z|\ge 3a\ge y_n.
\]
Taking the infimum over $z\in\partial\Omega$,
we obtain
\[
\delta(\Phi_i(y))
\ge\frac{y_n}{\sqrt{1+L^2}}.
\]
Since $(y',\varphi_i(y'))\in\partial\Omega$,
we also have
$\delta(\Phi_i(y))\le y_n$. Therefore,
\begin{equation}\label{distance-flatten-Nov}
\frac{y_n}{\sqrt{1+L^2}}
\le\delta(\Phi_i(y))\le y_n,
\qquad y\in D_{r_i}^+.
\end{equation}
Moreover, $|\det D\Phi_i|=1$ almost
everywhere, and
\[
|\nabla(u\circ\Phi_i)(y)|
\le C|\nabla u(\Phi_i(y))|.
\]
Since $\mu(x)=\delta(x)^\alpha$ on $U_i$,
it follows from \eqref{distance-flatten-Nov}
that
\begin{equation}\label{weight-flatten-Nov}
cy_n^\alpha\le\mu(\Phi_i(y))
\le Cy_n^\alpha,
\qquad y\in D_{r_i}^+.
\end{equation}
Here, $c>0$ and $C>0$ depend only on
$\alpha$ and the Lipschitz character
of $\Omega$.

\smallskip
We first establish the weighted
Sobolev-Poincar\'e inequality on each $U_i$.
For a boundary patch $U_i$, let us denote
\[
w_i(y)=u(\Phi_i(y)),
\quad\text{and}\quad
c_i=
\left(\int_{D_{r_i}^+}y_n^\alpha\,dy\right)^{-1}
\int_{D_{r_i}^+}w_i(y)y_n^\alpha\,dy.
\]
Applying the weighted Sobolev-Poincar\'e inequality
(Lemma \ref{weighted-sobolev-poincare}) on $D_{r_i}^+$ and using
\eqref{weight-flatten-Nov}, we obtain
\[
\begin{aligned}
\left(\int_{U_i}|u(x)-c_i|^q\,d\mu(x)\right)^{2/q}
&\le C\left(\int_{D_{r_i}^+}
|w_i(y)-c_i|^qy_n^\alpha\,dy\right)^{2/q}\\
&\le Cr_i^{2-N(1-2/q)}
\int_{D_{r_i}^+}|\nabla w_i(y)|^2
y_n^\alpha\,dy\\
&\le Cr_i^{2-N(1-2/q)}
\int_{U_i}|\nabla u(x)|^2\,d\mu(x).
\end{aligned}
\]
Similarly, applying the weighted Poincar\'e
inequality on $D_{r_i}^+$, we have
\[
\int_{U_i}|u(x)-c_i|^2\,d\mu(x)
\le Cr_i^2\int_{U_i}|\nabla u(x)|^2\,d\mu(x).
\]

\smallskip
Now, let $U_i=B_{r_i}(x_i)$ be an interior
ball. From \eqref{bound-mu}, and since $\delta(x)>\bar{\rho}$ on $U_i$,
it follows from the definition of $\mu$
and the assumption on the transition region
that
\[
0<c_0 \le\mu(x)\le C_0<\infty,
\qquad x\in U_i.
\]
Since $Q\ge n$, the exponent $q$ belongs
to the range of the usual Sobolev-Poincar\'e
inequality on $U_i$. Then, applying the
Sobolev-Poincar\'e and Poincar\'e inequalities
on $U_i$, we obtain
\[
\left(\int_{U_i}|u(x)-c_i|^q\,d\mu(x)\right)^{2/q}
\le C\int_{U_i}|\nabla u(x)|^2\,d\mu(x),
\]
and
\[
\int_{U_i}|u(x)-c_i|^2\,d\mu(x)
\le C\int_{U_i}|\nabla u(x)|^2\,d\mu(x),
\]
where $c_i$ denotes the integral average
of $u$ over $U_i$. 

\smallskip
By the estimates on the interior and the boundary patches $\{U_i\}_{i=1}^m$ that we just proved, and
$2-N(1-2/q)\ge0$, we conclude that
there is $C>0$ such that
\begin{equation}\label{local-SP-weight-mu-Nov}
\left(\int_{U_i}|u(x)-c_i|^q\,d\mu(x)\right)^{2/q}
+\int_{U_i}|u(x)-c_i|^2\,d\mu(x)
\le C\int_{U_i}|\nabla u(x)|^2\,d\mu(x)
\end{equation}
for every $i=1, 2, \ldots,m$.

\smallskip
We next establish the weighted Poincar\'e
inequality on $\Omega$. Since $\Omega$ is
connected, the finite covering
$\{U_i\}_{i=1}^m$ can be connected by
chains of overlapping patches.
Let $U_i$ and $U_j$ satisfy
$U_i\cap U_j\ne\emptyset$. Since
$\mu(U_i\cap U_j)>0$, we have
\[
\begin{aligned}
\mu(U_i\cap U_j)|c_i-c_j|^2
&\le2\int_{U_i\cap U_j}|u(x)-c_i|^2\,d\mu(x)
+2\int_{U_i\cap U_j}|u(x)-c_j|^2\,d\mu(x)\\
&\le C\int_\Omega|\nabla u(x)|^2\,d\mu(x),
\end{aligned}
\]
where we used \eqref{local-SP-weight-mu-Nov}
in the last estimate.
Since the number of patches is finite,
by connecting each $U_i$ to a fixed patch
$U_1$, we obtain
\begin{equation}\label{patch-means-Nov}
|c_i-c_1|^2
\le C\int_\Omega|\nabla u(x)|^2\,d\mu(x),
\qquad i=1,\ldots,m.
\end{equation}
Then, applying \eqref{local-SP-weight-mu-Nov}
and \eqref{patch-means-Nov}, we infer that
\[
\begin{aligned}
\int_\Omega|u(x)-c_1|^2\,d\mu(x)
&\le\sum_{i=1}^m\int_{U_i}|u(x)-c_1|^2\,d\mu(x)\\
&\le2\sum_{i=1}^m\int_{U_i}|u(x)-c_i|^2\,d\mu(x)
+2\sum_{i=1}^m\mu(U_i)|c_i-c_1|^2\\
&\le C\int_\Omega|\nabla u(x)|^2\,d\mu(x).
\end{aligned}
\]
Since $\bar u$ is the weighted mean of $u$,
we obtain
\begin{equation}\label{global-P-weight-mu-Nov}
\int_\Omega|u(x)-\bar u|^2\,d\mu(x)
\le\int_\Omega|u(x)-c_1|^2\,d\mu(x)
\le C\int_\Omega|\nabla u(x)|^2\,d\mu(x).
\end{equation}

\smallskip
Now, let $\{\eta_i\}_{i=1}^m$ be a smooth
partition of unity subordinate to the
covering $\{U_i\}_{i=1}^m$, satisfying
\[
0\le\eta_i(x)\le1,
\qquad
\operatorname{supp}\eta_i\subset U_i,
\qquad
\sum_{i=1}^m\eta_i(x)=1
\quad\text{in }\Omega.
\]
Using \eqref{local-SP-weight-mu-Nov} and
\eqref{global-P-weight-mu-Nov}, we have
\[
\begin{aligned}
\mu(U_i)|c_i-\bar u|^2
&\le2\int_{U_i}|u(x)-c_i|^2\,d\mu(x)
+2\int_{U_i}|u(x)-\bar u|^2\,d\mu(x)\\
&\le C\int_\Omega|\nabla u(x)|^2\,d\mu(x).
\end{aligned}
\]
Since $\mu(U_i)>0$ and the covering is finite,
we infer that
\begin{equation}\label{means-control-mu-Nov}
|c_i-\bar u|
\le C\left(\int_\Omega|\nabla u(x)|^2\,d\mu(x)
\right)^{1/2},
\qquad i=1,\ldots,m.
\end{equation}
Note that
\[
u(x)-\bar u
=\sum_{i=1}^m\eta_i(x)(u(x)-c_i)
+\sum_{i=1}^m\eta_i(x)(c_i-\bar u).
\]
Then, applying the triangle inequality,
\eqref{local-SP-weight-mu-Nov}, and
\eqref{means-control-mu-Nov}, we obtain
\[
\begin{aligned}
\left(\int_\Omega|u(x)-\bar u|^q\,d\mu(x)\right)^{1/q}
&\le\sum_{i=1}^m
\left(\int_{U_i}|\eta_i(x)(u(x)-c_i)|^q\,d\mu(x)
\right)^{1/q}\\
&\quad+\sum_{i=1}^m|c_i-\bar u|
\left(\int_{U_i}|\eta_i(x)|^q\,d\mu(x)\right)^{1/q}\\
&\le C\left(\int_\Omega|\nabla u(x)|^2\,d\mu(x)
\right)^{1/2}.
\end{aligned}
\]
Consequently,
\begin{equation}\label{global-SP-weight-mu-unit}
\left(\int_\Omega|u(x)-\bar u|^q\mu(x)\,dx
\right)^{2/q}
\le C\int_\Omega|\nabla u(x)|^2\mu(x)\,dx
\end{equation}
when $\operatorname{diam}(\Omega)=1$.

\smallskip
Finally, let $R=\operatorname{diam}(\Omega)$
and make the change of variables
$x=x_0+Ry$. Let us denote
\[
\widehat\Omega=\frac{\Omega-x_0}{R},
\quad
\widehat u(y)=u(x_0+Ry),
\quad\text{and}\quad
\widehat\mu(y)=R^{-\alpha}\mu(x_0+Ry).
\]
Note that $\operatorname{diam}(\widehat\Omega)=1$.
Moreover, if
\[
\widehat\delta(y)
=\operatorname{dist}(y,\partial\widehat\Omega),
\]
then
\[
\widehat\mu(y)=\widehat\delta(y)^\alpha
\quad\text{if}\quad
\widehat\delta(y)\le r_0/R,
\]
while $\widehat\mu$ is bounded above and
below by positive constants away from the
boundary. Thus, the preceding argument
applies to $\widehat\Omega$ with the weight
$\widehat\mu$. Since $N=n+\alpha$, we have
\[
\begin{aligned}
\int_\Omega|u(x)-\bar u|^q\mu(x)\,dx
&=R^N\int_{\widehat\Omega}
|\widehat u(y)-\overline{\widehat u}|^q
\widehat\mu(y)\,dy,\\
\int_\Omega|\nabla u(x)|^2\mu(x)\,dx
&=R^{N-2}\int_{\widehat\Omega}
|\nabla\widehat u(y)|^2\widehat\mu(y)\,dy.
\end{aligned}
\]
Applying \eqref{global-SP-weight-mu-unit}
on $\widehat\Omega$, we obtain
\[
\left(\int_\Omega|u(x)-\bar u|^q\mu(x)\,dx
\right)^{2/q}
\le CR^{2-N(1-2/q)}
\int_\Omega|\nabla u(x)|^2\mu(x)\,dx.
\]
This completes the proof.
\end{proof}
\subsection{H\"{o}lder regularity estimates for homogeneous equations} \label{simple-eqn} 
 For some $x_0= (x_{0}',x_{0n}) \in \bR^{n-1}\times \bR$ and some fixed $\rho>0$, we denote the upper-half cylinder in $\bR^n$ by
\[D_{\rho}^{+}(x_0)=D_{\rho}(x_0)\cap \bR_{+}^{n},\]
where $D_\rho(x_0)=B_{\rho}'(x_{0}')\times(x_{0n}-\rho,x_{0n}+\rho)$ and $B_{\rho}'(x_{0}')$ denotes the ball in $\bR^{n-1}$ centered at $x_{0}'$ with radius $\rho$. 
We also denote the flat boundary part of $D_\rho^{+}(x_0)$ by
\[\Gamma_{\rho}(x_0)=\p D_{\rho}^{+}(x_0)\cap (\bR^{n-1}\times \{0\}).\]
Note that if $\Gamma_{\rho}(x_0) \neq \emptyset$, then $\Gamma_{\rho}(x_0) = {\overline{B_{\rho}'(x_{0}')} \times \{0\}}$. When $x_0 = 0$, we simply write $D_{\rho}^{+}(0) = D_{\rho}^{+}$, $D_{\rho}(0) = D_{\rho}$, and $\Gamma_{\rho}(0)=\Gamma_{\rho}$.

\smallskip
{Let $r>0$ and $x_0\in\overline{\bR_+^n}$.} Let $\bB: {D_r^{+}(x_0)}\times \bR^n \rightarrow \bR^n$ be a vector field. Assume that there exists a constant $\Lambda \geq 1$ satisfying
\begin{equation}\label{bB-growth}
	\Lambda^{-1}|\xi-\zeta|^2 \leq \langle \bB(x, \xi)-\bB(x,\zeta), \xi -\zeta \rangle,
	\quad \text{and} \quad
	|\bB(x,\xi)| \leq \Lambda |\xi|,
\end{equation}
for all $x \in {D_r^{+}(x_0)}$ and all $\xi, \zeta \in \bR^n$. Moreover, for $\alpha \in (-1,\infty)$, let $\omega: \bR^n  \rightarrow \bR_{+}$ 
be a weight defined by 
\begin{equation}\label{omega-first}
	\omega(x) = |x_{n}|^{\alpha}(1-h(x))^{\alpha}, \quad \forall\, x = (x',x_n)\in \bR^{n-1} \times \bR,
\end{equation}
where $h: \bR^n \rightarrow [0,1)$ is a measurable function satisfying 
$\|h\|_{L_\infty(\bR^n)}\leq\tfrac12$.  By Remark~\ref{weight-remark} and $\|h\|_\infty\le1/2$, we have $\omega\in A_{2+\max\{\alpha,0\}}$ and $[\omega]_{2+\max\{\alpha,0\}}\le C(n,\alpha)$.

\smallskip
We consider the following class of equations:
\begin{equation}\label{boundary-eqn}
	\left\{
	\begin{array}{rcll}
		\diver(\omega(x)\bB(x,\nabla g)) &=& 0 & \text{in } {D_r^{+}(x_0)},\\[0.6em]
		\langle \omega(x)\bB(x,\nabla g), \overrightarrow{e_n}\rangle& = &0 
		& \text{on } {\Gamma_r(x_0)},
	\end{array}
	\right.
\end{equation}

\begin{definition}
	Let $r>0$ {and $x_0\in\overline{\bR_+^n}$}. Let $\omega$ be defined as in \eqref{omega-first}. Assume that $\bB: {D_r^{+}(x_0)}\times \bR^n\rightarrow \bR^n$ satisfies \eqref{bB-growth}. We say that $g\in W^{1,2}({D_{r}^{+}(x_0)},\omega)$ is a \emph{weak solution} of \eqref{boundary-eqn} 
	if it satisfies
	\begin{equation*}
		\int_{{D_{r}^{+}(x_0)}}\omega(x)\langle \bB(x, \nabla g(x)), \nabla\varphi(x)\rangle\,dx = 0,
	\end{equation*}
	for all $\varphi\in C^{\infty}_0({D_{r}(x_0)})$.
\end{definition}

\smallskip
Next, throughout the paper, for some non-empty open set $B$ and for a given function $f \in L^1(B, \gamma)$ with some given locally integrable weight function $\gamma: B \rightarrow [0, \infty)$, we  denote its mean  in $B$ with respect to the weight $\omega$ by
\begin{equation} \label{evarage-def}
\fint_{B}f(x)d\gamma(x) = \dfrac{1}{\gamma(B)}\int_{B}f(x)d\gamma(x), 
\end{equation}
where
\[
d\gamma(x) = \gamma(x) dx \quad \text{and} \quad \gamma(B) =\int_{B} \gamma(x) dx \in (0, \infty).
\]
Applying the weighted Poincar\'e and Sobolev-Poincar\'e inequalities (Lemmas \ref{weighted-poincare-mean} and \ref{weighted-sobolev-poincare}), we can use the De Giorgi-Nash-Moser argument as in \cite{FKS} to obtain the following H\"older estimates. 
\begin{lemma}\label{holder3}
	Let $\Lambda\geq 1$, $\alpha\in(-1, \infty)$. Then there exist $C=C(n,\Lambda,\alpha)>0$ and  $\beta=\beta(n,\Lambda, \alpha) \in(0,1)$ such that the following holds. Let $\bB: {D_r^{+}(x_0)}\times \bR^n \rightarrow \bR^n$ satisfy \eqref{bB-growth} with some $r>0$ and $x_0\in\overline{\bR_+^n}$, and $\omega$ satisfies  \eqref{omega-first}. Assume that $g\in W^{1,2}({D_r^+(x_0)},\omega)$ is a weak solution of \eqref{boundary-eqn}. Then it follows that
	\[\|g\|_{L^{\infty}({D_{r/2}^{+}(x_0)})} \leq C \Big(\fint_{{D_r^{+}(x_0)}} |g(x)|^2 d\omega(x)\Big)^{1/2},\]
	and
	\[|g(x)-g(y)| \leq C  {\Big(\fint_{{D_r^+(x_0)}}|g(x)|^2\,d\omega\Big)^{1/2}} \Big(\dfrac{|x-y|}{r}\Big)^{\beta}, \quad\forall\ x,y \in \overline{{D_{r/2}^{+}(x_0)}}.\]
\end{lemma}
\noindent
Similarly, applying the weighted Poincar\'e and Sobolev-Poincar\'e inequalities (Lemmas \ref{weighted-poincare-mean} and \ref{weighted-sobolev-poincare}) and Gehring's lemma, we obtain the following Meyers-Gehring type self-improving estimate. We skip the proof as it is standard now. See \cite{ME,KK} for the uniformly elliptic case.
\begin{lemma}\label{meyers_reverse3}
	Let $\Lambda\geq 1$, $\alpha\in (-1,\infty)$. Then, there exists $\overline{\delta}=\overline{\delta}(\Lambda,n,\alpha)>0$ such that the following holds. Assume that $\bB: {D_r^{+}(x_0)}\times {\bR^n} \rightarrow \bR^n$ satisfies \eqref{bB-growth} with some $r>0$ and $x_0\in\overline{\bR_+^n}$, and $\omega$ satisfies  \eqref{omega-first}. If $g\in W^{1,2}({D_r^{+}(x_0)},\omega)$ is a weak solution of \eqref{boundary-eqn}, then for every $p\in [2,2+\overline{\delta}]$, there is $C=C(\Lambda, n, 
	p, \alpha)>0$ such that
	\[\Big(\fint_{{D_{r/2}^{+}(x_0)}} |\nabla g(x)|^{p} d\omega(x)\Big)^{\frac{1}{p}}\leq C \Big(\fint_{{D_{r}^{+}(x_0)}} |\nabla g(x)|^{2}d\omega(x)\Big)^{1/2}.\]
\end{lemma}

\subsection{Lipschitz estimates and stability of simple equations}\label{subsec:homogeneous-lipschitz}\label{subsec:lipschitz-approximation}
We recall the following result from \cite[Lemma 3.3]{Dong-Phan-Sire}.
For $r>0$ and $x_0=(x_0',x_{0n})\in\overline{\bR_+^n}$, let
\[
\bA_0:((x_{0n}-r)_+,x_{0n}+r)\times\bR^n\longrightarrow\bR^n
\]
is measurable in $x_n$, Lipschitz continuous in $\xi$, and satisfies
\begin{equation} \label{A-zero-cond}
\left\{
\begin{split}
& \langle\bA_0(x_n,\xi)-\bA_0(x_n,\zeta),\xi-\zeta\rangle
\geq\Lambda^{-1}|\xi-\zeta|^2, \quad \forall \ \xi,\zeta\in\bR^n,\\ 
& \bA_0(x_n,0) =0,\quad
|D_{\xi}\bA_0(x_n,\xi)|\leq\Lambda, \quad\forall \ \xi \in\bR^n\setminus \{0\}
\end{split} \right.
\end{equation}
for all {$x_n \in ((x_{0n}-r)_+,x_{0n}+r)$}, where $s_+ = \max\{0, s\}$.
\begin{lemma}\label{homogeneous-lipschitz}
Let $\Lambda\geq1$, $\alpha>-1$, $r>0$, and $x_0\in\overline{\bR_+^n}$. Suppose that \eqref{A-zero-cond} holds and $u\in W^{1,2}(D_r^+(x_0),x_n^\alpha)$ satisfies
\[
\int_{D_r^+(x_0)}
\langle\bA_0(x_n,\nabla u(x)),\nabla\varphi(x)\rangle x_n^\alpha\,dx=0,
\qquad \forall \ \varphi\in C_0^\infty(D_r(x_0)),
\]
then
\[
\|\nabla u\|_{L^\infty(D_{r/2}^+(x_0))}
\leq C\left(
\frac{1}{\int_{D_r^+(x_0)}x_n^\alpha\,dx}
\int_{D_r^+(x_0)}|\nabla u(x)|^2x_n^\alpha\,dx
\right)^{1/2},
\]
where $C=C(\Lambda,n,\alpha)>0$.
\end{lemma}

\smallskip
Next, we recall the following stability result from \cite[Proposition 4.3]{Dong-Phan-Sire}. To do so, we need {some notation}. {For $\rho>0$ and $x_0\in\overline{\bR_+^n}$, suppose that $\bB:D_{2\rho}^+(x_0)\times\bR^n\to\bR^n$} is measurable in $x$, Lipschitz continuous in $\xi$, and satisfies
\begin{equation} \label{B-cond-simple}
\left\{
\begin{split}
& \langle\bB(x,\xi)-\bB(x,\zeta),\xi-\zeta\rangle
\geq\Lambda^{-1}|\xi-\zeta|^2, \quad \forall\ \xi,\zeta\in\bR^n,\\
& \bB(x,0) =0,\qquad |D_\xi\bB(x,\xi)|\leq\Lambda \quad \forall \  \xi \in\bR^n \setminus\{0\},
\end{split} \right.
\end{equation}
for all {$x \in D_{2\rho}^+(x_0)$}. In the spirit of \cite{MR2338417}, let us denote the weighted partial mean oscillation of $\bB$ in $D_\rho^+(x_0)$ by
\begin{equation} \label{theta-hat.def}
\widehat\Theta_{\rho,x_0}(\bB,\omega)
=\fint_{D_\rho^+(x_0)}\sup_{\xi\in\bR^n\setminus\{0\}} \frac{|\bB(x,\xi)-\overline{\bB}_{B_\rho'(x_0')}(x_n,\xi)|}{|\xi|}\,d\omega(x),
\end{equation}
where
\[
\overline{\bB}_{B_\rho'(x_0')}(x_n,\xi)
=\fint_{B_\rho'(x_0')}\bB(x',x_n,\xi)\,dx'.
\]
Here, $\omega(x)=|x_n|^\alpha\bar\omega(x)$  for $x = (x', x_n) \in \bR^{n-1}\times \bR$ and $\bar\omega : \bR^n \rightarrow \bR$ is a measurable, non-negative bounded function, and $B_{\rho}'(x_0') = \{y' \in \bR^{n-1}: {|y' -x_0'|} < \rho\}$.
\begin{proposition}\label{lipschitz-approximation-DPS}
Let $\Lambda\geq1$, $\alpha>-1$, $\rho>0$, and $x_0\in\overline{\bR_+^n}$. There exist sufficiently small numbers
$\bar{\delta}_0=\bar{\delta}_0(\Lambda,n,\alpha)>0$ and
$\sigma_0=\sigma_0(\Lambda,n,\alpha)>0$ such that the following holds. Assume that \eqref{B-cond-simple} holds on $D_{2\rho}^+(x_0)$,  $\|1-\bar\omega\|_{L^\infty(D_{2\rho}^+(x_0))}\leq \bar{\delta}_0$, and $v\in W^{1,2}(D_{2\rho}^+(x_0),\omega)$ and $\bG\in L^2(D_{2\rho}^+(x_0),\omega)^n$ satisfy
\[
\int_{D_{2\rho}^+(x_0)}\langle\bB(x,\nabla v(x)),\nabla\varphi(x)\rangle\,d\omega(x)
=\int_{D_{2\rho}^+(x_0)}\langle\bG(x),\nabla\varphi(x)\rangle\,d\omega(x),
\qquad \forall \ \varphi\in C_0^\infty(D_{2\rho}(x_0)).
\]
Then there exists
$w\in W^{1,2}(D_\rho^+(x_0),\omega)$ such that, for every $\sigma\in(0,\sigma_0)$,
\[
\begin{aligned}
&\left(\fint_{D_\rho^+(x_0)}|\nabla v(x)-\nabla w(x)|^2\,d\omega(x)\right)^{1/2}\\
&\quad\leq C\left(\widehat\Theta_{\rho,x_0}(\bB,\omega)^{\frac{\sigma}{2(2+\sigma)}}
+\|1-\bar\omega\|_{L^\infty(D_{2\rho}^+(x_0))}\right)
\left(\fint_{D_{2\rho}^+(x_0)}|\nabla v(x)|^2\,d\omega(x)\right)^{1/2}\\
&\qquad+C\left(\fint_{D_{2\rho}^+(x_0)}|\bG(x)|^2\,d\omega(x)\right)^{1/2},
\end{aligned}
\]
and
\[
\begin{aligned}
\|\nabla w\|_{L^\infty(D_{\rho/2}^+(x_0))}
&\leq C\left(\fint_{D_{2\rho}^+(x_0)}|\nabla v(x)|^2\,d\omega(x)\right)^{1/2} +C\left(\fint_{D_{2\rho}^+(x_0)}|\bG(x)|^2\,d\omega(x)\right)^{1/2},
\end{aligned}
\]
where $C=C(\Lambda,\alpha,n,\sigma)>0$.
\end{proposition}

\section{Boundary gradient estimates in flat domains}\label{sec:Lp-flat}

	As in Section \ref{simple-eqn}, for some $x_0= (x_{0}',x_{0n}) \in \bR^{n-1}\times \bR$ and some fixed $\rho>0$, we denote the upper-half cylinder in $\bR^n$ by
	\[D_{\rho}^{+}(x_0)=D_{\rho}(x_0)\cap \bR_{+}^{n},\]
	where $D_\rho(x_0)=B_{\rho}'(x_{0}')\times(x_{0n}-\rho,x_{0n}+\rho)$ and $B_{\rho}'(x_{0}')$ denotes the ball in $\bR^{n-1}$ centered at $x_{0}'$ with radius $\rho$. 
	We also denote the flat boundary part of $D_\rho^{+}(x_0)$ by
\[\Gamma_{\rho}(x_0)=\p D_{\rho}^{+}(x_0)\cap (\bR^{n-1}\times \{0\}).\]
	Note that if $\Gamma_{\rho}(x_0) \neq \emptyset$, then $\Gamma_{\rho}(x_0) = {\overline{B_{\rho}'(x_{0}')} \times \{0\}}$. When $x_0 = 0$, we simply write $D_{\rho}^{+}(0) = D_{\rho}^{+}$, $D_{\rho}(0) = D_{\rho}$, and $\Gamma_{\rho}(0)=\Gamma_{\rho}$.
	
\smallskip
Let $\bK\subset \bR$ be a given open interval, and let $\bA = (\A_1,\A_2,...,\A_n)$ be a vector field
 \[ 
 \begin{array}{cccc}
 \bA:  & D_{2}^{+}\times{\overline{\bK}}\times \bR^n & \rightarrow  & \bR^n \\
 \ & (x, z, \xi) & \mapsto & \bA(x,z, \xi).
 \end{array}
\] 
Assume that $\bA$ is measurable with respect to $x\in D_{2}^{+}$, Lipschitz in $z\in\mathbb{K}$, and  assume that there is $\Lambda \geq 1$ so that 

\begin{equation}\label{A-growth}
	\left\{\begin{array}{ll}
 |\xi-\zeta|^2 \leq \Lambda \langle \bA(x,z,\xi)-\bA(x,z,\zeta), \xi-\zeta \rangle &\quad \forall \  z\in \overline{\bK}, \ \xi,\zeta\in\bR^n\\ 
		|\bA(x,z_1,\xi)-\bA(x,z_2,\xi)|\leq \Lambda  |z_1-z_2| |\xi| &\quad \forall \  z_1,z_2\in \overline{\bK}, \ \xi \in\bR^n\\
\bA(x,z,0)=0, \ |D_{\xi}\bA(x,z,\xi)|\leq\Lambda
        & \quad \forall \  z\in\overline{\bK},\ \xi \in\bR^n \setminus\{0\}
        \end{array} \right.
\end{equation}
for a.e.  $x\in D_{2}^{+}$. For $\alpha\in(-1,\infty)$, let $\omega: \bR^n \rightarrow \bR_{+}$ be a weight defined by 
\begin{equation}\label{omega}
	\omega(x) = |x_{n}|^{\alpha}(1-h(x))^{\alpha}, \quad \forall \ x = (x',x_n)\in \bR^{n-1} \times \bR,
\end{equation}
where $h: \bR^n \rightarrow [0,1)$ is a measurable function satisfying $\|h\|_{L^\infty(\bR^n)}\leq\tfrac12$.

\smallskip
For every $\lambda > 0$, we study the following class of equations 
\begin{equation}\label{flat-eqs}
\left\{
\begin{array}{cccl}
		\diver(\omega(x)\bA(x,\lambda u,\nabla u)) &= &\diver(\omega(x)\bF(x)) &\text{in } D_{2}^{+}\\
	\lim\limits_{x\rightarrow x_0}	\omega(x)\Big(\A_{n}(x,\lambda u,\nabla u)-\F_n(x)\Big) &= &0 & x_0 \in \Gamma_{2}.
	\end{array} \right.
\end{equation}
where $\bF=(\F_1,\F_2,...,\F_n): D_{2}^{+} \rightarrow \bR^n$ is a given measurable vector field. For a given number $p\in(1,\infty)$, we say that $u\in W^{1,p}(D_{2}^{+},\omega)$ is a weak solution of \eqref{flat-eqs} if $\lambda u (x) \in \mathbb{K}$ for a.e. $x \in D_2^+$, and
\begin{equation*}
	\int_{D_{2}^{+}}\omega(x)\langle \bA(x,\lambda u(x),\nabla u(x)),  \nabla\varphi(x)\rangle dx = \int_{D_{2}^{+}}\omega(x)\langle \bF(x), \nabla\varphi(x)\rangle dx,
\end{equation*}
for any $\varphi\in C^{\infty}_0(D_{2})$.  As in \eqref{theta-hat.def}, we denote the weighted partial mean oscillation in $D_\rho^+(x_0)$ by
\begin{equation} \label{theta-hat-A.def}
\widehat\Theta_{\rho,x_0}(\bA,\omega)
=\fint_{D_\rho^+(x_0)}\sup_{z\in\bK}\sup_{\xi\in\bR^n\setminus\{0\}}
\frac{|\bA(x,z, \xi)-\overline{\bA}_{B_\rho'(x_0')}(x_n, z, \xi)|}{|\xi|}\,d\omega(x).
\end{equation}

\smallskip
The main result of this section is the following theorem, which provides a local boundary $L^p$ estimate for gradients of weak solutions of \eqref{flat-eqs}.
\begin{theorem}\label{flat-thrm}
	Let $\Lambda\geq 1$, $M>0$, $p\in[2,\infty)$, and $\alpha\in(-1, \infty)$. There exists a sufficiently small positive number $\overline{\delta}_0 = \overline{\delta}_0(\Lambda, M,\alpha,p,n)$ such that the following assertion holds. Assume that $\bA:D_2^{+}\times\overline{\bK}\times {\bR^n} \rightarrow \bR^n$  satisfies \eqref{A-growth} and $\omega$ satisfies \eqref{omega}, $\|h\|_{L^\infty(\bR^n)}\leq \overline{\delta}_0$, and
	\begin{equation}\label{A-bmo-flat}
		\sup_{0<\rho<R_0}\sup_{x\in\overline{ D_{1}^{+}}}\widehat \Theta_{\rho,x}(\bA,\omega)\leq \overline{\delta}_0,
	\end{equation}
	for some $R_0\in (0,1)$. Then, if $u\in W^{1,2}(D_{2}^{+},\omega)$ is a weak solution of \eqref{flat-eqs} satisfying $[[\lambda u]]_{\textup{BMO}_{R_0}({D_2^+},\omega)}\leq M$ with $\bF\in L^p(D_{2}^{+},\omega)^n$ and $\lambda>0$, we have $\nabla u \in L^p(D_1^+, \omega)$ and
\begin{equation}\label{flat-estimate}
	\|\nabla u\|_{L^p(D_{1}^{+},\omega )}^{p}\leq C\Big( \omega(D_2^{+})^{1-p/2}\|\nabla u\|_{L^2(D_{2}^{+},\omega )}^{p} + \|\bF\|_{L^p(D_{2}^{+},\omega )}^{p} \Big),
\end{equation}
where $C=C(\Lambda,M,p,\alpha,n,R_0)>0$.
\end{theorem}

\smallskip
We use the level-set {argument} introduced in \cite{CP} to prove the theorem. As in \cite{CP} and the developments in \cite{Buyn-Wang2, Buyn-Wang, BPS, Dong-Phan-Sire, GN, NP, P, P2}, Theorem \ref{flat-thrm} follows from the Proposition \ref{prop-approx} below using the level-set argument and a Krylov-Safonov covering lemma.  Because the argument is now standard, we skip the proof and only focus on the proof of Proposition \ref{prop-approx}. We note that the statement of the Proposition \ref{prop-approx} below, all average integrals are understood in the sense as in \eqref{evarage-def}.
\begin{proposition}\label{prop-approx} Let $\Lambda\geq 1$, $M>0$,  $\alpha\in (-1, \infty)$. For every $\epsilon\in(0,1)$, there exist small numbers $\kappa(\Lambda, M, n, \alpha, \epsilon)>0$ and $\delta_0(\Lambda, M, n, \alpha, \epsilon)\in(0,\epsilon)$ such that the following holds. Assume that $\bA:D_2^+\times\overline{\bK}\times\bR^n\rightarrow\bR^n$ satisfies \eqref{A-growth} and $\omega$ satisfies \eqref{omega} with $\|h\|_{L^\infty(\bR^n)}\leq\delta_0$, 
and
\begin{equation}
	\sup_{\rho\in(0,1)}\sup_{x\in \overline{D_{1}^{+}}}\widehat\Theta_{\rho,x}(\bA,\omega)<\delta_0 \quad \text{and} \quad  \fint_{D_{4r}^{+}(x_0)}|\bF(x)|^2d\omega(x)\leq \delta_0^2
\end{equation}
for some $x_0\in \overline{D_{1}^{+}}$ and some $r\in(0,1/4)$. Assume also that $u\in W^{1,2}(D_{2}^{+},\omega)$ is a weak solution of \eqref{flat-eqs} with some $\lambda>0$  and $u$ satisfies $[[\lambda u]]_{\textup{BMO}_{1}({D_2^+},\omega)}\leq M$, 
\begin{align*}
 \fint_{D_{4r}^{+}(x_0)}|\nabla u(x)|^2d\omega(x) \leq 1,  \quad \text{and} \quad  \fint_{D_{8\kappa r}^{+}(x_0)}|\nabla u (x)|^2d\omega(x) \leq 1.
	\end{align*}
Then, there is a function $w\in W^{1,2}(D_{2\kappa r}^{+}(x_0),\omega)$ such that the following estimates hold
	
	\begin{equation}\label{u-u3-bound}
		\Big(\fint_{D_{2\kappa r}^{+}(x_0)}|\nabla u(x) - \nabla w(x)|^2d\omega(x)\Big)^{1/2}\leq \epsilon,
	\end{equation}
and \begin{equation}\label{u3-bound}
	\|\nabla w\|_{L^{\infty}(D_{\kappa r}^{+}(x_0))}\leq  C(\Lambda, \alpha, n).
\end{equation}
\end{proposition}

\smallskip
The rest of this section is to prove Proposition \ref{prop-approx}. We divide the proof of Proposition \ref{prop-approx} into two comparison steps and then apply Proposition \ref{lipschitz-approximation-DPS}. In the first step, we freeze the dependence of $\bA$ on the given solution $u$, replace $\omega$ by $x_n^\alpha$, and approximate the weak solution $u$ of \eqref{flat-eqs} by the solution $u_1$ of the homogeneous equation \eqref{u1-eq} below. We note as in \eqref{evarage-def}, below
\[
\bar{u}_B = \fint_{B} u(x) d\omega(x) = \frac{1}{\omega(B)} \int_{B} u(x) \omega(x)\, dx
\]
for any measurable set $B$ with $\omega(B) \in (0, \infty)$.
\begin{lemma}\label{app-lemma1} Let $\Lambda\geq 1$, $M>0$, $\alpha\in(-1,\infty)$, and fix $\kappa\in(0,1)$. For every $\epsilon\in(0,1)$, there exists $\delta_1=\delta_1(\Lambda, n, \alpha, \epsilon,M)>0$ sufficiently small such that the following holds. Assume that $\bA:D_2^+\times\overline{\bK}\times\bR^n\rightarrow\bR^n$ and $\omega$ satisfy \eqref{A-growth} and \eqref{omega}, respectively, and that
	\begin{equation*}
\|h\|_{L^\infty(\bR^n)}\leq \delta_1 \quad \text{and}\quad 	\fint_{D_{r}^{+}(x_0)} |\bF(x)|^2 d\omega(x)\leq \delta_1,
	\end{equation*}
	for some ${x_0\in\overline{D_{1}^{+}}}$, $r\in(0,1)$.
	Then, for every $\lambda>0$, if $u\in W^{1,2}(D_{2}^{+},\omega)$  is a weak solution of \eqref{flat-eqs} satisfying
	
	\begin{equation*}
	\fint_{D_{r}^{+}(x_0)} |\nabla u(x)|^2 d\omega(x)\leq 1 \quad \text{and} \quad  \lambda\Big(\fint_{D_{r}^{+}(x_0)} | u-\bar{u}_{D_{r}^{+}(x_0)} |^2d\omega(x)\Big)^{1/2}\leq M,
	\end{equation*}
then
	\begin{equation}\label{l1-s1}
	\fint_{D_{r}^{+}(x_0)}|\nabla u(x) -\nabla u_1(x)|^2d\omega(x) \leq \epsilon^2,
	\end{equation}
where $u_1 \in W^{1,2}({D_{r}^{+}(x_0)},\omega)$ is the weak solution of
	\begin{equation}\label{u1-eq}
	\left\{	\begin{array}{cccl}
			\diver(x_n^\alpha \bA(x,\lambda u,\nabla u_1)) &=& 0 &\text{ in } D_{r}^{+}(x_0)\\
			u_1 &= & u - \bar{u}_{D_{\kappa r}^{+}(x_0)} &\text{ on } \p D_{r}^{+}(x_0) \setminus \Gamma_{r}(x_0)\\
		\lim\limits_{x_n\rightarrow 0^+}	x_{n}^{\alpha}\bA_{n}(x,\lambda u,\nabla u_1) &= &0  &\text{ on } \Gamma_{r}(x_0).
		\end{array}	\right.		
	\end{equation}
Moreover, it holds that
	\begin{equation}\label{l1-s2}
		\lambda\left(\fint_{D_{r}^{+}(x_0)}  \big|u_1(x)-\bar u_{D_r^+(x_0)}+\bar u_{D_{\kappa r}^+(x_0)}\big|^2 d\omega(x)\right)^{1/2}\leq C(n,\alpha)(\epsilon r \lambda + M)
	\end{equation}
	and
	\begin{equation}\label{l1-s3}
		\Big(\fint_{D_{r}^{+}(x_0)} |\nabla u_1(x)|^2d\omega(x)\Big)^{1/2}\leq 2.
	\end{equation}
\end{lemma}
\begin{proof} We note that by applying the Minty-Browder theorem (see \cite[Theorem 1.5]{MB} and \cite[Lemma 4.4]{Dong-Phan-Sire}), we obtain the existence and uniqueness of $u_1$. Testing \eqref{flat-eqs} with $\varphi=u_1-(u-\bar u_{D_{\kappa r}^+(x_0)})$, we obtain
		\begin{equation*}
		\int_{D_{r}^{+}(x_0)}\langle \bA(x,\lambda u, \nabla u(x)), \nabla u(x) - \nabla u_1(x) \rangle \omega(x) dx = \int_{D_{r}^{+}(x_0)}\langle\bF,\nabla u(x) - \nabla u_1(x)  \rangle \omega(x) dx.
	\end{equation*}
Similarly, testing \eqref{u1-eq} with $\varphi=u_1-(u-\bar u_{D_{\kappa r}^+(x_0)})$, we obtain
\begin{equation*}
\int_{D_{r}^{+}(x_0)}\langle \bA(x,\lambda u, \nabla u_1(x)), \nabla u(x) - \nabla u_1(x) \rangle x_n^{\alpha} dx = 0.
\end{equation*}
Therefore, from \eqref{A-growth}, H\"{o}lder's inequality and Young's inequality,  we obtain
\begin{align*}
		& \int_{D_{r}^{+}(x_0)} |\nabla u(x) - \nabla u_1(x)|^2 \omega(x) dx \\
		&\leq \Lambda \int_{D_{r}^{+}(x_0)}\langle \bA(x,\lambda u, \nabla u(x))-\bA(x,\lambda u, \nabla u_1(x)), \nabla u(x) - \nabla u_1(x) \rangle \omega(x) dx\\
		& = \Lambda \int_{D_{r}^{+}(x_0)}\langle\bF,\nabla u(x) - \nabla u_1(x)  \rangle \omega(x) dx \\
		&\quad - \Lambda \int_{D_{r}^{+}(x_0)} \langle \bA(x,\lambda u, \nabla u_1(x)), \nabla u(x) - \nabla u_1(x) \rangle (\omega(x)-x_n^{\alpha}) dx\\
		&\leq C(\Lambda)  \int_{D_{r}^{+}(x_0)}|\bF(x)|^2\omega(x) dx + \dfrac{1}{3} \int_{D_{r}^{+}(x_0)} |\nabla u(x) - \nabla u_1(x)|^2 \omega(x) dx \\
		& \quad + C(\Lambda) \sup_{x\in D_{r}^{+}(x_0)}\Big| 1- \frac{x_n^{\alpha}}{\omega(x)}\Big|^2 \int_{D_{r}^{+}(x_0)}|\nabla u_1(x)|^2\omega(x) dx,
	\end{align*}
for $C = C(\Lambda, n)>0$. Then,
\begin{align*}
\int_{D_{r}^{+}(x_0)} |\nabla u(x) - \nabla u_1(x)|^2 \omega(x) dx & \leq  C(\Lambda) \int_{D_{r}^{+}(x_0)}|\bF(x)|^2\omega(x) dx \\
		&\quad + C(\Lambda) \sup_{x\in D_{r}^{+}(x_0)}\Big| 1- \frac{x_n^{\alpha}}{\omega(x)}\Big|^2 \int_{D_{r}^{+}(x_0)}|\nabla u_1(x)|^2\omega(x) dx.
	\end{align*}
Observe that from \eqref{omega} and as $\|h\|_{L^\infty}$ is sufficiently small, we have
\[
\Big| 1- \frac{x_n^{\alpha}}{\omega(x)}\Big| = |1 - (1- h(x))^{-\alpha}| \leq C \|h\|_{L^\infty(\bR^n)}.
\]
with $C = C(\alpha, n)>0$. Hence,
\begin{align*}
		&\int_{D_{r}^{+}(x_0)} |\nabla u(x) - \nabla u_1(x)|^2 \omega(x) dx\\
		& \leq  C(\Lambda)  \int_{D_{r}^{+}(x_0)}|\bF(x)|^2\omega(x) dx + C(\Lambda,\alpha)\|h\|_{L^\infty(\bR^n)}^{2} \int_{D_{r}^{+}(x_0)} |\nabla u_1(x)|^2\omega(x) dx\\
		&\leq  C(\Lambda)  \int_{D_{r}^{+}(x_0)}|\bF(x)|^2\omega(x) dx + C(\Lambda,\alpha)\|h\|_{L^\infty(\bR^n)}^{2} \int_{D_{r}^{+}(x_0)} |\nabla u(x)|^2\omega(x) dx\\
		& \quad + C(\Lambda,\alpha)\|h\|_{L^\infty(\bR^n)}^{2} \int_{D_{r}^{+}(x_0)} |\nabla u(x)-\nabla u_1(x)|^2\omega(x) dx\\
		&\leq  C(\Lambda)  \int_{D_{r}^{+}(x_0)}|\bF(x)|^2\omega(x) dx + C(\Lambda,\alpha)\|h\|_{L^\infty(\bR^n)}^{2} \int_{D_{r}^{+}(x_0)} |\nabla u(x)|^2\omega(x) dx\\
		&\quad 	+ \dfrac{1}{3} \int_{D_{r}^{+}(x_0)} |\nabla u(x)-\nabla u_1(x)|^2\omega(x) dx,
	\end{align*}
where we used $\delta_1>0$ sufficiently small so that $ C(\Lambda,\alpha)\|h\|_{L^\infty(\bR^n)}^{2} \leq C\delta_1^2 \leq 1/3$. Therefore, it follows that
	\begin{align*}
		\int_{D_{r}^{+}(x_0)} |\nabla u(x) - \nabla u_1(x)|^2 \omega(x) dx &\leq  C(\Lambda)  \int_{D_{r}^{+}(x_0)}|\bF(x)|^2\omega(x) dx\\
		&\quad + C(\Lambda,\alpha)\|h\|_{L^\infty(\bR^n)}^{2} \int_{D_{r}^{+}(x_0)} |\nabla u(x)|^2\omega(x) dx.
	\end{align*}
Next, we choose $\delta_1$ smaller, if necessary, so that $C(\Lambda)\delta_1 < \epsilon^2/2$ and $C(\Lambda,\alpha)\delta_1^2 < \epsilon^2/2$. Then, from the previous estimate and the choice of $\delta_1$, we infer that
	\begin{align*}
		&\dfrac{1}{\omega(D_r^{+}(x_0))}\int_{D_{r}^{+}(x_0)} |\nabla u(x) - \nabla u_1(x)|^2 d\omega(x)\\
		& \leq \dfrac{C(\Lambda)}{\omega(D_r^{+}(x_0))}\int_{D_{r}^{+}(x_0)}|\bF(x)|^2 d\omega(x)  + \dfrac{C(\Lambda,\alpha)\|h\|_{L^\infty(\bR^n)}^{2} }{\omega(D_r^{+}(x_0))}\int_{D_{r}^{+}(x_0)} |\nabla u(x)|^2 d\omega(x) \\
		& \leq C(\Lambda)\delta_1 + C(\Lambda,\alpha)\delta_1^2\\
		& \leq \epsilon^2.
	\end{align*}
This proves \eqref{l1-s1}.  For \eqref{l1-s2}, the boundary condition gives $u_1-(u-\bar u_{D_{\kappa r}^+(x_0)})=0$ on the nonflat boundary. By the triangle inequality, the Poincar\'e inequality in Lemma~\ref{weighted-poincare-top}, and the comparison of $d\omega$ with $x_n^\alpha\,dx$, we obtain
\begin{align*}
 &\lambda\Big(\fint_{D_r^+(x_0)}
 \big|u_1(x)-\bar u_{D_r^+(x_0)}+\bar u_{D_{\kappa r}^+(x_0)}\big|^2d\omega\Big)^{1/2}\\
 &\leq\lambda\Big(\fint_{D_r^+(x_0)}
 |u_1(x)-(u(x)-\bar u_{D_{\kappa r}^+(x_0)})|^2d\omega\Big)^{1/2}
 +\lambda\Big(\fint_{D_r^+(x_0)}|u(x)-\bar u_{D_r^+(x_0)}|^2d\omega\Big)^{1/2}\\
 &\leq C(n,\alpha)\lambda r
 \Big(\fint_{D_r^+(x_0)}|\nabla u_1(x)-\nabla u(x)|^2d\omega\Big)^{1/2}+M
 \leq C(n,\alpha)(\epsilon r\lambda+M).
\end{align*}
Finally, to prove \eqref{l1-s3}, we apply the triangle inequality and the assumptions of the lemma to obtain
	\begin{align*}
		&\Big( \fint_{D_{r}^{+}(x_0)} | \nabla u_1(x)|^2 d\omega(x)\Big)^{1/2}\\
		&\leq \Big(\fint_{D_{r}^{+}(x_0)} | \nabla u(x)|^2 d\omega(x)\Big)^{1/2} + \Big( \fint_{D_{r}^{+}(x_0)} | \nabla u(x) - \nabla u_1(x)|^2 d\omega(x)\Big)^{1/2}\\
		&\leq 1 + \epsilon \leq 2.
	\end{align*}
This proves \eqref{l1-s3} and completes the proof.
\end{proof}

\smallskip
Our next lemma is the most technical one in the paper.


\begin{lemma}\label{app-lemma2}
	Let $\Lambda\geq 1$, $M>0$, $\alpha\in(-1,\infty)$. For every $\epsilon\in(0,1)$, there exist sufficiently small positive numbers $\kappa = \kappa(\Lambda, M, \alpha, n, \epsilon)$ and  $\delta_2=\delta_2(\Lambda, M, \alpha,  n, \epsilon)$ such that the following holds. 
	Assume that $\bA$ satisfies \eqref{A-growth} and $\omega$ satisfies \eqref{omega}. Also, assume that 
	\begin{equation}
		\|h\|_{L^\infty(\bR^n)}\leq\delta_2 \quad\text{and}\quad \fint_{D_{r}^{+}(x_0)} |\bF(x)|^2 d\omega(x)\leq \delta_2,
	\end{equation}
	for some $x_0\in \overline{D}_{1}^{+}$, and  $r\in(0,1)$. Then, for every $\lambda>0$, if $u\in W^{1,2}(D_2^{+},\omega)$  is a weak solution of \eqref{flat-eqs} satisfying
\begin{equation}
\begin{split}
& \fint_{D_{r}^{+}(x_0)} |\nabla u(x)|^2 d\omega(x)\leq 1, \quad
\fint_{D_{2\kappa r}^{+}(x_0)} |\nabla u(x)|^2 d\omega(x)\leq 1,
\end{split}
\end{equation}
and $[[\lambda u]]_{\textup{BMO}_{1}({D_2^+},\omega)}\leq M,$ then the following estimates hold:
\begin{equation}\label{l2-s1}
	\fint_{D_{\kappa r}^{+}(x_0)}|\nabla u(x) -\nabla u_2(x)|^2d\omega(x) \leq \epsilon^2,
	\end{equation}
and
	\begin{equation}\label{l2-s2}
		\Big(	\fint_{D_{\kappa r}^{+}(x_0)}|\nabla u_2(x)|^2d\omega(x)\Big)^{1/2} \leq C(\alpha,n),
	\end{equation}
where $u_2 \in W^{1,2}({D_{\kappa r}^{+}(x_0)},\omega)$ is the weak solution of
	\begin{equation}\label{u2-eq}
			\left\{\begin{array}{cccl}
			\diver(x_n^\alpha \bA(x,\lambda \bar{u}_{ D_{\kappa r}^{+}(x_0)},\nabla u_2)) &=& 0& \text{ in } D_{\kappa r}^{+}(x_0)\\
			u_2 &= &u_1  &\text{ on } \p D_{\kappa r}^{+}(x_0) \setminus \Gamma_{\kappa r}(x_0)\\
		\lim\limits_{x_n\rightarrow 0}	x_{n}^{\alpha}\A_{n}(x,\lambda \bar{u}_{ D_{\kappa r}^{+}(x_0)},\nabla u_2) &=& 0  &\text{ on } \Gamma_{\kappa r}(x_0)
		\end{array}\right.			
\end{equation}
in which $u_1\in W^{1,2}(D_r^+(x_0),\omega)$ is the solution of \eqref{u1-eq} supplied by Lemma \ref{app-lemma1} with the same $\kappa$.
\begin{proof}
Let $\epsilon>0$. Let $\epsilon'\in(0,\epsilon/2)$, and $\kappa >0$ sufficiently small which will be determined later. By Lemma \ref{lm-doubling}-(i), there are $C_d=C_d(n,\alpha)\ge1$ and $\vartheta=\vartheta(n,\alpha)>0$ such that
\[
\frac{\omega(D_r^+(x_0))}{\omega(D_\rho^+(x_0))}
\le C_d\left(\frac r\rho\right)^\vartheta,\qquad 0<\rho\le r.
\]
By the local John--Nirenberg inequality,
\[
\lambda\Big(\fint_{D_r^+(x_0)}|u-\bar u_{D_r^+(x_0)}|^2\,d\omega\Big)^{1/2}
\le C(n,\alpha)M.
\]
Let $\delta_2=\delta_1(\Lambda,n,\alpha,\epsilon'\kappa^{\vartheta/2}/\sqrt{C_d},C(n,\alpha)M)>0$, where $\delta_1$ is defined in Lemma~\ref{app-lemma1}.	Applying Lemma~\ref{app-lemma1} with the chosen $\kappa$, we obtain $u_1\in W^{1,2}(D_r^+(x_0),\omega)$ such that
		\begin{equation}\label{l2-p1}
			\Big(\fint_{D_{r}^{+}(x_0)}|\nabla u(x) -\nabla u_1(x)|^2d\omega(x)\Big)^{1/2} \leq {\epsilon'\kappa^{\vartheta/2}/\sqrt{C_d}},
		\end{equation}
		and
		\begin{equation}\label{l2-p2}
			\lambda	\Big(\fint_{D_{r}^{+}(x_0)}{\big|u_1(x)-\bar u_{D_r^+(x_0)}+\bar u_{D_{\kappa r}^+(x_0)}\big|^2}d\omega(x)\Big)^{1/2} \leq C(n,\alpha)[r\epsilon'\kappa^{\vartheta/2}\lambda+M].
		\end{equation}
Moreover, 
		$$\Big(\fint_{D_{r}^{+}(x_0)}|\nabla u_1(x)|^2d\omega(x)\Big)^{1/2}\leq 2.$$
		Notice that if $x_{0n}>r$, then $D_{r}^{+}(x_0)=D_r(x_0)$. By \eqref{l2-p1} and triangle inequality, we get
		\begin{align*}
			&\Big(\fint_{D_{2\kappa r}^{+}(x_0)}|\nabla u_1(x)|^2 d\omega(x)\Big)^{1/2}\\
			& \leq \Big(\fint_{D_{2\kappa r}^{+}(x_0)}|\nabla u(x)|^2 d\omega(x)\Big)^{1/2} +\Big(\fint_{D_{2\kappa r}^{+}(x_0)}|\nabla u(x) - \nabla u_1(x)|^2 d\omega(x)\Big)^{1/2}\\
			& \leq 1 + \Big({\dfrac{C_d}{(2\kappa)^{\vartheta}}}\fint_{D_{ r}^{+}(x_0)}|\nabla u(x) - \nabla u_1(x)|^2 d\omega(x)\Big)^{1/2}\\
			&\leq 1 + \dfrac{\epsilon'}{2^{\vartheta/2}} \leq 2.
		\end{align*}
	{Applying the Caccioppoli estimate with the constant $\bar u_{D_r^+(x_0)}-\bar u_{D_{\kappa r}^+(x_0)}$ subtracted, and using \eqref{l2-p2}, we obtain}
		\begin{align*}
			\Big(\fint_{D_{2\kappa r}^{+}(x_0)}|\nabla u_1(x)|^2 d\omega(x) \Big)^{1/2}&\leq \dfrac{C(\Lambda,\alpha,n)}{(1-2\kappa)\kappa^{\vartheta/2}r}\Big(\fint_{D_{r}^{+}(x_0)}{\big|u_1-\bar u_{D_r^+(x_0)}+\bar u_{D_{\kappa r}^+(x_0)}\big|^2}d\omega(x)\Big)^{1/2}\\
			& \leq C(\Lambda,\alpha,n) \dfrac{1}{(1-2\kappa)\kappa^{\vartheta/2}r} \Big(r\epsilon'\kappa^{	\vartheta/2}+\dfrac{M}{\lambda}\Big)\\
			&\leq C(\Lambda,\alpha,n)\Big(\epsilon'+\dfrac{M}{\lambda \kappa^{\vartheta/2}r}\Big).
		\end{align*}
{The functions $u_1$ and $u_2$ solve \eqref{u1-eq} and \eqref{u2-eq}, respectively.} Thus, we can use $u_1-u_2$ as a test function which implies
		\begin{align*}
			&\int_{D_{\kappa r }^{+}(x_0)} \langle \bA(x,\lambda u(x), \nabla u_1(x)), \nabla u_2(x) - \nabla u_1(x)\rangle x_n^{\alpha} dx\\
			& = \int_{D_{\kappa r }^{+}(x_0)} \langle \bA(x,\lambda \bar{u}_{D_{\kappa r}^{+}(x_0)}, \nabla u_2(x)), \nabla u_2(x) - \nabla u_1(x)\rangle x_n^{\alpha} dx\\
			& =0.
		\end{align*}
{Since $\|h\|_\infty\leq1/2$, $d\omega$ and $x_n^\alpha\,dx$ are comparable; we use this comparison in the estimates below.}
Applying the structural condition \eqref{A-growth}, we infer that
		\begin{align*}
			&\int_{D_{\kappa r }^{+}(x_0)} |\nabla u_1(x) - \nabla u_2(x)|^2d\omega(x)\\
			&\leq C \int_{D_{\kappa r }^{+}(x_0)} \langle \bA(x,\lambda  \bar{u}_{D_{\kappa r}^{+}(x_0)}, \nabla u_1(x) )-\bA(x,\lambda  \bar{u}_{D_{\kappa r}^{+}(x_0)}, \nabla u_2(x) ), \nabla u_1(x) - \nabla u_2(x)\rangle {x_n^\alpha\,dx}\\
			& = C \int_{D_{\kappa r }^{+}(x_0)} \langle \bA(x,\lambda  \bar{u}_{D_{\kappa r}^{+}(x_0)}, \nabla u_1(x) ), \nabla u_1(x) - \nabla u_2(x)\rangle {x_n^\alpha\,dx}\\
			&\leq C \int_{D_{\kappa r }^{+}(x_0)} |\nabla u_1(x)|.|\nabla u_1(x) - \nabla u_2(x)|d\omega(x)\\
			&\leq C \int_{D_{\kappa r }^{+}(x_0)}|\nabla u_1(x)|^2d\omega(x) + \dfrac{1}{2}\int_{D_{\kappa r }^{+}(x_0)}|\nabla u_1(x) -\nabla u_2(x)|^2d\omega(x)
\end{align*}
where $C = C(\Lambda,n,\alpha)>0$. Then, by canceling the last term, we deduce that
		$$\int_{D_{\kappa r }^{+}(x_0)} |\nabla u_1(x) - \nabla u_2(x)|^2d\omega(x)\leq {C(\Lambda,n,\alpha)}\int_{D_{\kappa r}^{+}(x_0)}|\nabla u_1(x)|^2d\omega(x).$$
This estimate and the previous estimate together result in 
		\begin{align*}
				\Big(\fint_{D_{\kappa r}^{+}(x_0)}|\nabla u_1(x) -\nabla u_2(x)|^2d\omega(x)\Big)^{1/2}
			& \leq C(\Lambda,\alpha,n)	\Big(\fint_{D_{2\kappa r}^{+}(x_0)}|\nabla u_1(x)|^2d\omega(x)\Big)^{1/2}\\
			& \leq C(\Lambda,\alpha,n)\Big(\epsilon'+\dfrac{M}{\lambda\kappa^{\vartheta/2}r}\Big).
		\end{align*}
Therefore, we are left to prove that the right hand side of the above estimate is small. Let us consider the following two cases.

\smallskip		
\noindent
\textbf{Case 1}: $\dfrac{ C(\Lambda,\alpha,n)M}{\lambda\kappa^{\vartheta/2}r} \leq \dfrac{\epsilon}{4}$. We choose $\epsilon'\in (0,\epsilon/2)$ such that
\[C(\Lambda, \alpha, n)(\epsilon')<\dfrac{\epsilon}{4},\]
then
\[	\Big(\fint_{D_{\kappa r}^{+}(x_0)}|\nabla u_1(x) -\nabla u_2(x)|^2d\omega(x)\Big)^{1/2} \leq \epsilon/2.\]
From this and \eqref{l2-p1}, and the triangle inequality, we get
		\begin{align*}
			&	\Big(\fint_{D_{\kappa r}^{+}(x_0)}|\nabla u(x) - \nabla u_2(x)|^2d\omega(x)\Big)^{1/2}\\
			& \leq 	\Big(\fint_{D_{\kappa r}^{+}(x_0)}|\nabla u(x) - \nabla u_1(x)|^2d\omega(x)\Big)^{1/2} +	\Big(\fint_{D_{\kappa r}^{+}(x_0)}|\nabla u_1(x) - \nabla u_2(x)|^2d\omega(x)\Big)^{1/2}\\
			& \leq \Big({\dfrac{C_d}{\kappa^{\vartheta}}}\fint_{D_{ r}^{+}(x_0)}|\nabla u(x) - \nabla u_1(x)|^2d\omega(x)\Big)^{1/2} +\epsilon/2\\
			&\leq \epsilon'+\epsilon/2 \leq \epsilon.
		\end{align*}
		
\smallskip
\noindent
\textbf{Case 2}: $\frac{ C(\Lambda, \alpha, n)M}{\kappa^{\vartheta/2}r\lambda} > \frac{\epsilon}{4}$. It is equivalent to $\kappa^{\vartheta/2}r\lambda \epsilon < 4C(\Lambda, \alpha, n)M$. Since $\epsilon'\in(0,\epsilon/2)$, we have  $\kappa^{\vartheta/2}r\lambda \epsilon' < C(\Lambda,n, \alpha, M)$.
		Then, by \eqref{l2-p2},
		\begin{align*}
& \lambda \Big( \fint_{D_{r}^{+}(x_0)}{\big|u_1(x)-\bar u_{D_r^+(x_0)}+\bar u_{D_{\kappa r}^+(x_0)}\big|^2}d\omega(x)\Big)^{1/2}  \leq {C(n,\alpha)}(r\epsilon'\kappa^{\vartheta/2}\lambda+M)\\
&\leq C(\Lambda, n, \alpha, M).
\end{align*}
Applying Holder's regularity theory (Lemma \ref{holder3}) to $\lambda(u_1-\bar u_{D_r^+(x_0)}+\bar u_{D_{\kappa r}^+(x_0)})$, whose oscillation equals that of $v=\lambda u_1$, there is $\beta\in(0,1)$ depending on ${\Lambda,\alpha,n}$ so that
\[
|v(x)-v(y)|\leq N(\Lambda, n, \alpha, M)\Big(\dfrac{|x-y|}{r}\Big)^{\beta}\quad \forall \  x,y \in \overline{D^{+}_{r/4}}(x_0).
\]
Denote $\hat{u}=u-\bar{u}_{D_{\kappa r}^{+}(x_0)}$. We employ the growth condition \eqref{A-growth} of the vector field $\bA$ and the weak solution formulation with an appropriate test function to infer
		\begin{align*}
			&\int_{D_{\kappa r}^{+}(x_0)}|\nabla u_1(x)-\nabla u_2(x)|^2d\omega(x)\\
			&\leq N \int_{D_{\kappa r}^{+}(x_0)} \langle \bA(x,\lambda \bar{u}_{D_{\kappa r}^{+}(x_0)}, \nabla u_1(x))-\bA(x,\lambda \bar{u}_{D_{\kappa r}^{+}(x_0)}, \nabla u_2(x)), \nabla u_1(x) -\nabla u_2(x) \rangle {x_n^\alpha\,dx}\\
			&= N \int_{D_{\kappa r}^{+}(x_0)} \langle \bA(x,\lambda \bar{u}_{D_{\kappa r}^{+}(x_0)}, \nabla u_1(x))-\bA(x,\lambda u, \nabla u_1(x)), \nabla u_1(x) -\nabla u_2(x) \rangle {x_n^\alpha\,dx}\\
			&\leq N\int_{D_{\kappa r }^{+}(x_0)} |\lambda \hat{u}|.|\nabla u_1(x)|.|\nabla u_1(x)-\nabla u_2(x)|d\omega(x)\\
			& \leq N\int_{D_{\kappa r}^{+}(x_0)}|\lambda \hat{u}|^2|\nabla u_1(x)|^2d\omega(x) +\dfrac{1}{2}\int_{D_{\kappa r}^{+}(x_0)}|\nabla u_1(x)-\nabla u_2(x)|^2d\omega(x),
		\end{align*}
		where we used Young's inequality to derive the last estimate, and $N = N(\Lambda, \alpha)>0$. Following this, by canceling the last term, we get
		\begin{align*}
			\fint_{D_{\kappa r}^{+}(x_0)}|\nabla u_1(x)-\nabla u_2(x)|^2d\omega(x) \leq {C(\Lambda,\alpha)}\fint_{D_{\kappa r}^{+}(x_0)}|\lambda \hat{u}|^2|\nabla u_1(x)|^2d\omega(x).
		\end{align*}
		For $q_1>2$, sufficiently close to $2$, we write $q_0=\dfrac{2q_1}{q_1-2}>2$. Then, using Holder's inequality and self-improving regularity estimates (Lemma \ref{meyers_reverse3}),
		we have
		\begin{align*}
			&\dfrac{1}{\omega(D_{\kappa r}^{+}(x_0))} \int_{D_{\kappa r}^{+}(x_0)} |\nabla u_1(x) -\nabla u_2(x)|^2d\omega(x)\\
			&\leq {C(\Lambda,\alpha)}\Big( \fint_{D_{\kappa r}^{+}(x_0)} |\lambda \hat{u}|^{q_0}d\omega(x) \Big)^{\frac{2}{q_0}} \Big( \fint_{D_{\kappa r}^{+}(x_0)} |\nabla u_1(x)|^{q_1}d\omega(x) \Big)^{\frac{2}{q_1}}\\
			&\leq {C(\Lambda,n,\alpha)}\Big( \fint_{D_{\kappa r}^{+}(x_0)} |\lambda \hat{u}|^{q_0}d\omega(x) \Big)^{\frac{2}{q_0}} 
		\Big(\fint_{D_{2\kappa r}^{+}(x_0)} |\nabla u_1(x)|^{2}d\omega(x) \Big).
		\end{align*}

\smallskip		
Now, we estimate the first term on the right hand side of the previous estimate. As noted after  Definition \ref{BMO}, John-Nirenberg lemma follows from  \cite[Theorem 5]{Muckenhoupt-Weeden-BMO} and  \cite[Theorem 4.1]{Buckley-JN}. Hence, 
\begin{align*}
\dfrac{1}{\omega(D_{\kappa r}^{+}(x_0))} \int_{D_{\kappa r}^{+}(x_0)} |\lambda \hat{u}|^{q_0}d\omega(x) & \leq \Big(  \fint_{D_{\kappa r}^{+}(x_0)} |\lambda \hat{u}|^{2}d\omega(x)\Big)^{1/2}\Big(  \fint_{D_{\kappa r}^{+}(x_0)} |\lambda \hat{u}|^{2(q_0-1)}d\omega(x)\Big)^{1/2}\\
			&\leq {C(n,q_0,\alpha)} [[\lambda u]]_{\textup{BMO}_{1}({D_2^+},\omega)}^{q_0-1} \Big( \fint_{D_{\kappa r}^{+}(x_0)} |\lambda \hat{u}|^{2}d\omega(x)\Big)^{1/2}\\
			&\leq {C(n,q_0,M,\alpha)}\Big(  \fint_{D_{\kappa r}^{+}(x_0)} |\lambda \hat{u}|^{2}d\omega(x)\Big)^{1/2}.
		\end{align*}
Therefore,
\begin{align*}
& \fint_{D_{\kappa r}^{+}(x_0)} |\nabla u_1(x) -\nabla u_2(x)|^2d\omega(x)\\
&\leq {C(\Lambda,n,\alpha)}\Big( \fint_{D_{\kappa r}^{+}(x_0)} |\lambda \hat{u}(x)|^{q_0}d\omega(x) \Big)^{\frac{2}{q_0}}
		 \Big( \fint_{D_{2\kappa r}^{+}(x_0)} |\nabla u_1(x)|^{2}d\omega(x) \Big)\\
		&\leq {C(\Lambda,n,q_0,M,\alpha)}\Big(  \fint_{D_{\kappa r}^{+}(x_0)} |\lambda \hat{u}(x)|^{2}d\omega(x) {\Big)^{1/q_0}}
		\Big(\fint_{D_{2\kappa r}^{+}(x_0)}|\nabla u_1(x)|^2d\omega(x)\Big).	
		\end{align*}
		Moreover, recall that $v=\lambda u_1$. With $\hat{u}=u-\bar{u}_{D_{\kappa r}^{+}(x_0)}$, we have
		\begin{align*}
& \fint_{D_{\kappa r}^{+}(x_0)}|\lambda\hat{u}(x)|^2d\omega(x) =  \fint_{D_{\kappa r}^{+}(x_0)}|\lambda\hat{u}(x) - \lambda \bar{\hat{u}}_{D_{\kappa r}^{+}(x_0)}|^2d\omega(x) \\
			& \leq C\Big(  \fint_{D_{\kappa r}^{+}(x_0)}|\lambda(\hat{u}(x)-u_1(x))|^2d\omega(x)\Big.
		 +  \fint_{D_{\kappa r}^{+}(x_0)}|\lambda(u_1(x)- \overline{u}_{1,D_{\kappa r}^{+}(x_0)})|^2d\omega(x)\\
			&\quad \Big. + \quad  \fint_{D_{\kappa r}^{+}(x_0)}
			|\lambda(\bar{\hat{u}}_{D_{\kappa r}^{+}(x_0)}-\bar{u}_{1,D_{\kappa r}^{+}(x_0)})|^2d\omega(x) \Big)\\
			&\leq C \Big( \dfrac{1}{\kappa^\vartheta}\fint_{D_{ r}^{+}(x_0)}|\lambda(\hat{u}(x)-u_1)|^2d\omega(x)
		+ \fint_{D_{\kappa r}^{+}(x_0)}|v(x)-{\bar{v}_{D_{\kappa r}^{+}(x_0)}}|^2d\omega(x)\Big).
		\end{align*}
With the chosen boundary value, $\hat u-u_1=0$ on $\partial D_r^+(x_0)\setminus\Gamma_r(x_0)$. Using Poincar\'e inequality in Lemma \ref{weighted-poincare-top} and the comparison of $d\omega$ with $x_n^\alpha\,dx$ for the first term {and the oscillation bound for $v$ for the second term} of the above estimate, we obtain
		\begin{align*}
			&	\Big( \fint_{D_{\kappa r}^{+}(x_0)} |\lambda\hat{u}(x)|^2d\omega(x)\Big)^{1/2}\\
			& \leq C\Big[\dfrac{\lambda^2r^2}{\kappa^\vartheta}\fint_{D_{ r}^{+}(x_0)}|\nabla \hat{u}-\nabla u_1(x)|^2d\omega(x) +\Big(\sup_{x,y\in D_{\kappa r}^{+}(x_0)}|v(x)-v(y)|\Big)^2 \Big]^{1/2}\\
			&\leq C \Big[\dfrac{\lambda r}{\kappa^{\vartheta/2}}\Big(\fint_{D_{ r}^{+}(x_0)}|\nabla u(x)-\nabla u_1(x)|^2d\omega(x)\Big)^{1/2} +\sup_{x,y\in D_{\kappa r}^{+}(x_0)}|v(x)-v(y)| \Big]\\ 
			&\leq {C(\Lambda,n,M,\alpha)}\Big[\dfrac{\lambda r}{\kappa^{\vartheta/2}}\epsilon'\kappa^{\vartheta/2} + \Big(\dfrac{\kappa r}{r}\Big)^{\beta}\Big]\\
			&\leq {C(\Lambda,n,M,\alpha)}(\lambda r \epsilon'+\kappa^{\beta}).
		\end{align*}
		From this, we can compute
		\begin{align*}
			\fint_{D_{\kappa r}^{+}(x_0)}|\nabla u_1(x)-\nabla u_2(x)|^2d\omega(x)
			&\leq  {C(n,q_0,M,\alpha)} \Big(\fint_{D_{\kappa r}^{+}(x_0)}|\lambda\hat{u}|^2d\omega(x)\Big)^{\frac{q_1-2}{2q_1}}\\
			&\leq {C(n,q_0,M,\Lambda,\alpha)}(\lambda r \epsilon'+\kappa^{\beta})^{\frac{q_1-2}{q_1}}\\
			&\leq {C(n,q_0,M,\Lambda,\alpha)}\Big(\dfrac{\epsilon'}{\epsilon\kappa^{\vartheta/2}}+\kappa^{\beta}\Big)^{\frac{q_1-2}{q_1}}.
		\end{align*}
		We can choose $\kappa$ small so that 
\[\kappa^{\beta}{< \tfrac12\Big[}\dfrac{\Big(\frac{\epsilon}{2}\Big)^2}{{C(n,q_0,M,\Lambda,\alpha)}}\Big]^{\frac{q_1}{q_1-2}}.\]
		In addition, we choose $\epsilon'$ small so that
\[\dfrac{\epsilon'}{\epsilon\kappa^{\vartheta/2}} {< \tfrac12\Big[}\dfrac{{\big(\frac{\epsilon}{2}\big)}^2}{{C(n,q_0,M,\Lambda,\alpha)}}\Big]^{\frac{q_1}{q_1-2}}.\]
		Thus,
	\[\fint_{D_{\kappa r}^{+}(x_0)}|\nabla u_1(x) -\nabla u_2(x)|^2d\omega(x) {<(\epsilon/2)^2}.\] {From this and \eqref{l2-p1}, as in Case 1, we get
\[
\Big(\fint_{D_{\kappa r}^{+}(x_0)}|\nabla u(x)-\nabla u_2(x)|^2d\omega(x)\Big)^{1/2}
\leq \epsilon'+\epsilon/2 \leq \epsilon.
\]
This proves \eqref{l2-s1}.}
		Next, we prove the last estimate in the lemma. Note that
		\begin{align*}
			&\Big( \fint_{D_{\kappa r}^{+}(x_0)}|\nabla u_2(x)|^2d\omega(x)\Big)^{1/2}\\
			&\leq \Big( \fint_{D_{\kappa r}^{+}(x_0)}|\nabla u(x)|^2d\omega(x)\Big)^{1/2} + \Big(\fint_{D_{\kappa r}^{+}(x_0)}|\nabla u(x) - \nabla u_2(x)|^2d\omega(x)\Big)^{1/2}\\
			&\leq  C(\alpha,n) \Big(2^\vartheta \fint_{D_{2\kappa r}^{+}(x_0)}|\nabla u(x)|^2d\omega(x)\Big)^{1/2} +\quad \epsilon\\
			&\leq C(\alpha,n)2^{\vartheta/2}+\epsilon\\
			&\leq C(\alpha,n),
		\end{align*}
where we used triangle inequality and assumptions given in this lemma.
\end{proof}
\end{lemma}


\smallskip
We are now ready to prove Proposition \ref{prop-approx}.

\begin{proof}[Proof of Proposition \ref{prop-approx}]
By Lemma \ref{lm-doubling}-\textup{(i)} and \eqref{l2-s2}, choose $\bar C_0=\bar C_0(\Lambda,n,\alpha)\geq1$ such that
\[
\left(\frac{\omega(D_{4\kappa r}^+(x_0))}
{\omega(D_{2\kappa r}^+(x_0))}\right)^{1/2}\leq \bar  C_0,
\qquad C(\alpha,n)\leq \bar C_0,
\]
where $C(\alpha,n)$ is the constant in \eqref{l2-s2}.
Applying Lemma \ref{app-lemma2} with $4r$ and $\epsilon/(2\bar C_0)$ in place of $r$ and $\epsilon$, respectively, we obtain $\kappa\in(0,1/8)$ and $u_2$ on $D_{4\kappa r}^+(x_0)$ such that
\[
\left(\fint_{D_{4\kappa r}^+(x_0)}
|\nabla u(x)-\nabla u_2(x)|^2\,d\omega(x)\right)^{1/2}
\leq\frac{\epsilon}{2\bar C_0},
\qquad
\left(\fint_{D_{4\kappa r}^+(x_0)}
|\nabla u_2(x)|^2\,d\omega(x)\right)^{1/2}\leq \bar C_0,
\]
provided $\delta_0\leq\min\{1,\delta_2(\Lambda,M,\alpha,n,\epsilon/(2\bar C_0)), \bar{\delta}_0\}$, where $\bar{\delta}_0$ is defined in Proposition \ref{lipschitz-approximation-DPS}.

\smallskip
Now, let us set
\[
\bB(x,\xi)=\bA\left(x,\lambda\bar u_{D_{4\kappa r}^+(x_0)},\xi\right), \quad \text{for} \quad x \in D_2^+, \quad \xi \in \bR^n.
\]
By \eqref{A-growth}, $\bB$ satisfies the hypotheses of Proposition \ref{lipschitz-approximation-DPS}. Moreover,
\[
\widehat\Theta_{2\kappa r,x_0}(\bB,\omega)\leq\delta_0.
\]
Applying Proposition \ref{lipschitz-approximation-DPS} with $\rho=2\kappa r$ to $u_2$ with weight $x_n^\alpha$, $\bF=0$, and $\sigma=\sigma_0/2$, where $\sigma_0$ is given in that proposition. Since $d\omega(x)\simeq x_n^\alpha\,dx$, we obtain $w\in W^{1,2}(D_{2\kappa r}^+(x_0),\omega)$ satisfying
\[
\left(\fint_{D_{2\kappa r}^+(x_0)}
|\nabla u_2(x)-\nabla w(x)|^2\,d\omega(x)\right)^{1/2}
\leq C(\Lambda,n,\alpha)\delta_0^{\frac{\sigma}{2(2+\sigma)}},
\]
and
\[
\|\nabla w\|_{L^\infty(D_{\kappa r}^+(x_0))}
\leq C(\Lambda,n,\alpha).
\]
Here we used
\[
\widehat\Theta_{2\kappa r,x_0}(\bB,x_n^\alpha)
\leq C(\alpha)\widehat\Theta_{2\kappa r,x_0}(\bB,\omega)
\leq C(\alpha)\delta_0.
\]
Choose $\delta_0\in(0,\epsilon)$ smaller, if necessary, so that
\[
C(\Lambda,n,\alpha)\delta_0^{\frac{\sigma}{2(2+\sigma)}}\leq\frac{\epsilon}{2}.
\] 
By {Lemma \ref{lm-doubling}-\textup{(i)}} and the triangle inequality,
\[
\begin{aligned}
&\left(\fint_{D_{2\kappa r}^+(x_0)}
|\nabla u(x)-\nabla w(x)|^2\,d\omega(x)\right)^{1/2}\\
&\quad\leq C_0\left(\fint_{D_{4\kappa r}^+(x_0)}
|\nabla u(x)-\nabla u_2(x)|^2\,d\omega(x)\right)^{1/2}
+\frac{\epsilon}{2}\leq\epsilon.
\end{aligned}
\]
This proves \eqref{u-u3-bound} and \eqref{u3-bound}.
\end{proof}
\section{Global gradient estimates and proof of Theorem \ref{main-thrm}}\label{sec:main-proofs} 
In this section, we prove the main theorem \ref{main-thrm}. There are three subsections which include interior, boundary, and global estimates.
\subsection{Interior estimates}\label{subsec:interior}

For every $r\in (0, \text{diam}(\Omega)/2)$, we define
\[\Omega^{r}=\{x\in \Omega: \text{dist}(x,\p \Omega)>r\}.\]
For $p\in[1,\infty)$, $r>0$, and $\lambda>0$, we say that $u\in W^{1,p}(\Omega^{r},\mu)$ is a weak solution of 
\begin{equation}\label{interior-eq}
	\diver(\mu(x)\bA(x,\lambda u, \nabla u )) = \diver(\mu(x)\bF(x))\quad\quad\text{in } \Omega^r
\end{equation}
if $\lambda u(x) \in {\overline{\bK}}$ and
\[\int_{\Omega^{r}} \langle \bA(x,\lambda u(x), \nabla u(x) ), \nabla \varphi\rangle \mu(x)dx  =\int_{\Omega^{r}} \langle \bF(x), \nabla \varphi(x)\rangle \mu(x)dx,\]
for any $\varphi \in C_0^{\infty}(\Omega^{r})$.  

\smallskip
Recall also the constants $\rho, c_0, C_0$ are defined in \eqref{bound-mu}. We are now ready to state the result on the interior domain. We want to note that this result is proved in \cite[Theorem 1.2]{P}.
\begin{theorem}\label{interior-thrm} For every $\Lambda\geq 1$, $M>0$, $p\in[2,\infty)$ and $\alpha\in(-1,\infty)$, there exists a sufficiently small positive number $\tau_1 = \tau_1(\Lambda,M,p,\alpha, n)$ such that the following assertion holds. Assume that $\bA:\Omega \times\overline{\bK}\times {\bR^n}\rightarrow \bR^n$  satisfies \eqref{bA-growth} and $\mu$ satisfies \eqref{mu-weight}, and
	\begin{equation}\label{A-bmo-interior}
		\sup_{0<\rho<R_0}\sup_{x\in\overline{\Omega}} \Theta_{\rho,x}(\bA,\mu)\leq \tau_1,
	\end{equation}
	for some $R_0\in (0,1)$. Then, for each $r\in [\rho,r_0)$, if $u\in W^{1,2}(\Omega^{r},\mu)$ is a weak solution of \eqref{interior-eq} satisfying 	$[[\lambda u]]_{\textup{BMO}_{R_0}(\Omega^{r},\mu)}\leq M$  with $\bF\in L^p(\Omega^{r},\mu)^n$, we have
	\begin{equation}\label{flat-estimate-interior}
		\|\nabla u\|_{L^p(\Omega^{2r},\mu )}^{p}\leq C\mu(\Omega^{r})^{1-p/2}\|\nabla u\|_{L^2(\Omega^{r},\mu )}^{p} + C\|\bF\|_{L^p(\Omega^{r},\mu )}^{p}
	\end{equation}
	for some constant $C=C(\Lambda, {R_0,r,r_0,M,p,\alpha,n,\Omega,c_0,C_0})>0$ {depending also on the modulus of continuity of $\mu$ on $\overline{\Omega^{\rho}}$}.
\end{theorem}
\begin{proof} The case $p=2$ is immediate. Hence, we only consider $p>2$. We apply \cite[Theorem 1.2]{P}. We then need to set up the problem in the setting in \cite[Theorem 1.2]{P} and check the ellipticity condition, the smallness condition on the mean oscillations of the coefficients.  Since $\mu$ is positive and continuous, for each $\varepsilon\in(0,1/3)$, we can choose $s\in(0,\min\{r/4,R_0/4\})$ such that
\[
\left|\frac{\mu(x)}{\mu(x_0)}-1\right|\leq\varepsilon
\quad\text{for }x_0\in\overline{\Omega^{2r}},
\quad x\in B_{2s}(x_0).
\]
Set $\widetilde\bA=\mu\bA/\mu(x_0)$ and
$\widetilde\bF=\mu\bF/\mu(x_0)$. Then $\widetilde\bA$ satisfies
\eqref{bA-growth} with $2\Lambda$ in place of $\Lambda$, and
\[
\diver\widetilde\bA(x,\lambda u,\nabla u)
=\diver\widetilde\bF
\quad\text{in }B_{2s}(x_0).
\]
Denote $\gamma(x) = \mu(x)/\mu(x_0)$. For $y \in B_s(x_0)$ and ${B=B_{r'}(y)}$ with $0<r' \leq s$, by a direct manipulation, we obtain
\[
\widetilde{\mathbf A}-\fint_B\widetilde{\mathbf A}\,dx
=\gamma (\mathbf A-\overline{\mathbf A}_B)
 +(\gamma-\gamma_B)\overline{\mathbf A}_B, \quad \text{where} \quad \gamma_B = \fint_{B} \gamma(x) dx.
\]
From this, it follows that
\begin{align*}
\Theta_{r',y}(\widetilde\bA,1) & \leq \gamma_{B} \Theta_{r', y}(\bA, \mu) + \Lambda \fint_{B} |\gamma(x) - \gamma_B|\, dx \\
& \leq 2\Theta_{r',y}(\bA,\mu)+2\Lambda\varepsilon \leq 2 \tau_1 + 2\Lambda \varepsilon.
\end{align*}
On the other hand, we also have
\begin{align*}
&\lambda\fint_B\Big|u-\fint_B u\,dx\Big|\,dx
\leq 2\lambda\fint_B|u-\bar u_B|\,dx
\leq 2\frac{\sup_B\mu}{\inf_B\mu}\lambda\fint_B|u-\bar u_B|\,d\mu
\leq \frac{2(1+\varepsilon)}{1-\varepsilon}M\leq 4M.
\end{align*}
We choose $\tau_1$ and $\varepsilon$ so small such that
$2\tau_1+2\Lambda\varepsilon$ is less than the smallness constant in
\cite[Theorem 1.2]{P} with the structural ellipticity constant $2\Lambda$ and BMO
bound $4M$. Applying that theorem, using
\[
0<c_0 \leq\mu(x)\leq C_0 <\infty
\quad\text{for }x\in\Omega^r,
\]
and covering $\overline{\Omega^{2r}}$ by finitely many balls
$B_s(x_0)$, we obtain \eqref{flat-estimate-interior}.
{By \eqref{mu-weight} and \eqref{bound-mu}, these bounds depend only on $r,r_0,\alpha,c_0,C_0$. The choice of $s$ also depends on the modulus of continuity of $\mu$ on $\overline{\Omega^{\rho}}$.}
\end{proof}

\subsection{Local boundary estimates}\label{subsec:boundary} 

Let $r_0>0$ be the small number defined in \eqref{mu-weight}. Since $\Omega$ is $(\delta,\rho_0)$-Lipschitz, it follows that for $R\in (0, \min(\rho_0,r_0)/2)$, and $x_0\in \p\Omega$, by Definition \ref{lipschitz-domain} and with a rotation and translation, 
\begin{equation}\label{flatten-C}
	C_{2R}(x_0) =\{x=(x',x_n)\in B_{2R}'(x'_0)\times\bR: \gamma(x')<x_n<\gamma(x') + 2R\}\subset \Omega
\end{equation}

\begin{equation}\label{flatten-Gamma}
	\Gamma_{2R}(x_0)= \{(x',\gamma(x')):x'\in B_{2R}'(x'_0) \}\subset \p \Omega
\end{equation}
where $\gamma: \overline{B_{2R}'}(x'_0)\rightarrow \bR$ is a Lipschitz function which satisfies
\[\gamma(0)=0,\quad \nabla \gamma(0)=0,\quad \text{and }\|\nabla \gamma\|_{L^{\infty}(B_{2R}'(x'_0))}\leq \delta.\]
Recall also that
\[\mu(x)=\text{dist}(x,\p\Omega)^{\alpha},\quad x\in C_{2R}(x_0).\]

\smallskip
We study the equation \eqref{main-eqs} in the neighborhood of $0\in\p\Omega$:
\begin{equation}\label{boundary-eqs}
	\left\{\begin{array}{cccl}
	\diver(\mu(x)\bA(x,\lambda u,\nabla u)) &=&\diver(\mu(x)\bF(x)) & \text{in } C_{2R}(x_0)\\
\lim\limits_{{x\rightarrow z}}	\mu(x)\Big(\bA(x,\lambda u,\nabla u) - \bF(x)\Big) \cdot \overrightarrow{\nu} &=& 0 & \text{when }{z\in \Gamma}_{2R}(x_0).
\end{array}\right.
\end{equation}
For $p\in(1,\infty)$, we say that $u\in W^{1,p}(C_{2R}(x_0),\mu)$ is a weak solution of \eqref{boundary-eqs}{ if $\lambda u(x)\in\overline{\bK}$ almost everywhere and}
\begin{equation}
\int_{C_{2R}(x_0)} \langle \bA(x,\lambda u(x),\nabla u(x)),\nabla \varphi(x)\rangle \mu(x)dx = \int_{C_{2R}(x_0)} \langle \bF(x),\nabla \varphi(x)\rangle \mu(x)dx
\end{equation}
for all $\varphi \in C^{\infty}(\overline{C_{2R}}(x_0))$ vanishing locally on $\p C_{2R}(x_0)\setminus \Gamma_{2R}(x_0)$.

\smallskip
We next establish a local boundary estimate for weak solutions of \eqref{boundary-eqs}. By translation, we take $x_0=0$ and write  $C_{2R}(0) = C_{2R}$, $ \Gamma_{2R}(0)= \Gamma_{2R}$, and $B_{2R}'(0)=B_{2R}'$.

\begin{theorem}\label{boundary-thrm} Let $R\in (0, \min\{\rho_0/8,r_0/2\})$. Let $M>0$. For each $p\in[2,\infty)$, there exists a sufficiently small constant $\tau_2=\tau_2(\Lambda,M,p,n,\alpha)>0$ such that the following holds. If growth condition \eqref{A-growth} holds, $\Omega$ is $(\delta,\rho_0)$-Lipschitz with $\delta\in(0,\tau_2)$ and 
\[{\Theta_{\rho,x}}(\bA,\mu)<\tau_2,\quad \forall \  x\in \overline{C_R}, \forall \  \rho\in(0,R_0)\]
for some $R_0\in(0,1)$, then for any weak solution $u\in W^{1,2}(C_{2R},\mu)$ of \eqref{boundary-eqs} for some {$\lambda>0$} satisfying 	$[[\lambda u]]_{\textup{BMO}_{R_0}({C_{2R}},\mu)}\leq M$, and $\bF\in L^p(C_{2R},\mu)^n$, we have $\nabla u \in L^p(C_R,\mu)$ and
\begin{equation}\label{boundary-result}
	\|\nabla u\|_{L^p(C_R,\mu)}^{p} \leq C\mu(C_{2R})^{1-p/2}\|\nabla u\|_{L^2(C_{2R},\mu)}^{p} + C	\|\bF\|_{L^p(C_{2R},\mu)}^{p},
\end{equation}
where $C=C(\Lambda,M,p,n,\alpha,\widehat R_0)>0$ with $\widehat R_0=\min\{R_0/(2R),1/4\}$.
\end{theorem}

\smallskip
In order to prove this theorem, we begin by flattening the boundary $\Gamma_{2R}$ and then apply Theorem \ref{flat-thrm}.

\smallskip
Let $\Phi: C_{2R} \rightarrow D_{2R}^{+}$ and $\Psi: D_{2R}^{+}\rightarrow C_{2R}$ be defined by
\begin{align*}
	\Phi(x)=(x',x_n-\gamma(x'))\quad \forall \  x = (x',x_n) \in C_{2R},\\
	\Psi(y)=(y',y_n+\gamma(y'))\quad \forall \  y = (y',y_n) \in D_{2R}^{+}.
\end{align*}
We observe that
\begin{equation}
	\nabla \Phi(x)=\left(\begin{array}{cc}
	I_{n-1} & 0\\
	-\p_{x'}\gamma(x') & 1
	\end{array}\right)
\end{equation}
{for almost all} $x=(x',x_n)\in C_{2R}$, and
\begin{equation}
		\nabla \Psi(y)=\left(\begin{array}{cc}
		I_{n-1} & 0\\
		\p_{y'}\gamma(y') & 1
	\end{array}\right)
\end{equation}
{for almost all} $y=(y',y_n)\in D_{2R}^{+}$, where $I_{n-1}$ is the $(n-1)\times(n-1)$ identity matrix. This implies $\text{det}(\nabla \Psi)=\text{det}(\nabla \Phi)=1$ and
\[\Phi=\Psi^{-1},\quad \nabla\Psi(y)=[\nabla\Phi(\Psi(y))]^{-1},\quad{\text{for a.e. }y}\in D_{2R}^{+}.\]
Moreover, as $\delta\in(0,1)$,
\begin{equation}\label{grads-psi-phi}
\begin{array}{l}
\|\nabla\Phi\|_{L^{\infty}(C_{2R})}^{2}\leq n+ \|\nabla\gamma\|_{L^{\infty}({B_{2R}}')}^{2}\leq n+\delta^2\leq n+1\\
\|\nabla\Psi\|_{L^{\infty}(D_{2R}^{+})}^{2}\leq n+1.
\end{array}
\end{equation}
Recall that 
\[T_{2R}={\overline{B_{2R}'}\times\{0\}}=\p D_{2R}^{+} \cap \{y_n=0\},\]
and denote
\begin{equation}\label{newA}
	\omega(y)=\mu(\Psi(y)) \quad\text{and}\quad \hat{\bA}(y,z,\xi)=\bA(\Psi(y),z,\xi[\nabla\Phi(\Psi(y))])[\nabla \Phi(\Psi(y))]^*
\end{equation}
for $y\in D_{2R}^{+}$. 

\smallskip
We are now ready to introduce the following lemma which is obtained via an appropriate change of variables.
\begin{lemma}\label{change-of-variables-lm}
	Assume that $\Omega$ is $(\delta,\rho_0)$-Lipschitz. If  $u\in W^{1,p}(C_{2R},\mu)$ is a weak solution of \eqref{boundary-eqs} and $w(y)=u(\Psi(y))$ for any $y\in D_{2R}^{+}$, then there is $w\in W^{1,p}(D_{2R}^{+},\omega)$ that is a weak solution of 
	\begin{equation}\label{scaled-eqs}
		\left\{\begin{array}{cccl}
			\diver[\omega(y)\mathrm{\hat{\bA}}(y,\lambda w,\nabla w)] &= & \diver(\omega(y)\bG(y)) \quad &\text{in } D_{2R}^{+},\\
		\lim\limits_{y_n\rightarrow 0}	\omega(y)\big(\hat{\A}_n(y,\lambda w,\nabla w)-G_n(y)\big) &=&0 \quad &\text{on }  T_{2R},
		\end{array}\right.
	\end{equation}
where $\bG(y)= \bF(\Psi(y))[\nabla\Phi(\Psi(y))]^{*}$ with $y\in D_{2R}^{+}$.
\end{lemma}

\begin{proof}
By the change of variables and as $u\in W^{1,p}(C_{2R},\mu)$, we obtain $w\in W^{1,p}(D_{2R}^+,\omega)$. Applying the same change of variables to the weak formulation for $u$, we obtain the weak formulation for $w$.
\end{proof}

\smallskip
Let us denote
\[[[\bA]]_{\textup{BMO}_{R_0}}(C_R,\mu) = \sup_{0<\rho<R_0}\sup_{x\in C_R}\Theta_{\rho,x}(\bA,\mu),\]
where $\Theta_{\rho,x}(\bA,\mu)$ is defined in {Definition \ref{A-bmo-def}}. A similar definition can be made for $\hat{\bA}$. The next result gives the estimate of the mean oscillation of $\hat{\bA}$ with the weight $\mu$. Its proof can be found in \cite{Dong-Phan-Sire}.


\begin{lemma}\label{A-comparision} Suppose that $\Omega$ is a $(\delta,\rho_0)$-Lipschitz domain and $\bA$ satisfies \eqref{A-growth}. For $0<\delta\leq1$, $\hat{\bA}$ satisfies the revised \eqref{A-growth} on $D_{2R}^{+}$ with $\Lambda$ replaced by $\Lambda_*=(n+1)\Lambda$. Moreover, there is $N_0=N_0(\Lambda,n)>0$ such that
\[[[\hat{\bA}]]_{\textup{BMO}_{R_0}}(D_{R}^{+},\omega)\leq N_0\big([[\bA]]_{\textup{BMO}_{2R_0}}(C_{R},\mu) + \delta\big),\]
where $\hat{\bA}$ is defined in \eqref{newA}.
\end{lemma}

\smallskip
To verify the structural assertion, put $J=\nabla\Phi(\Psi(y))$. By \eqref{grads-psi-phi}, $\|J\|^2,\|J^{-1}\|^2\leq n+1$. The row-vector formula \eqref{newA} gives $\hat\bA(y,z,0)=0$ and
\[
 |\hat\bA(y,z,\xi)-\hat\bA(y,z,\zeta)|
 \leq\Lambda\|J\|^2|\xi-\zeta|\leq\Lambda_*|\xi-\zeta|,
\]
\[
 |\hat\bA(y,z_1,\xi)-\hat\bA(y,z_2,\xi)|
 \leq\Lambda_*|z_1-z_2|\,|\xi|.
\]
The growth bound follows by setting $\zeta=0$. Moreover,
\begin{align*}
 &\langle\hat\bA(y,z,\xi)-\hat\bA(y,z,\zeta),\xi-\zeta\rangle\\
 &=\langle\bA(\Psi(y),z,\xi J)-\bA(\Psi(y),z,\zeta J),(\xi-\zeta)J\rangle\\
 &\geq\Lambda^{-1}|(\xi-\zeta)J|^2
 \geq\Lambda_*^{-1}|\xi-\zeta|^2.
\end{align*}
Finally, since $\|J-I\|\leq\delta$,
\[
 |\bA(\Psi(y),z,\xi J)J^*-\bA(\Psi(y),z,\xi)|
 \leq\Lambda(\|J\|+1)\delta|\xi|.
\]
This supplies the structural $O(\delta)$ error in the coefficient-oscillation comparison. Hence constants introduced by this coordinate change depend only on $n$ and the enlarged $\Lambda$.

\smallskip
We are ready to give the proof for Theorem \ref{boundary-thrm}.
\begin{proof}[Proof of Theorem \ref{boundary-thrm}.]
{We take $R<\min\{\rho_0/8,r_0/2\}$ and $\delta\leq1$. For $x=\Psi(y)\in C_{2R}$, any boundary point outside the graph chart is farther from $x$ than $y_n$. Since $\|\nabla\gamma\|_\infty\leq\delta$, we have
\[
\frac{y_n}{\sqrt{1+\delta^2}}
\leq\operatorname{dist}(\Psi(y),\partial\Omega)\leq y_n.
\]
Thus, with $\widetilde h(y)=1-\operatorname{dist}(\Psi(y),\partial\Omega)/y_n$,
\[
\omega(y)=y_n^\alpha(1-\widetilde h(y))^\alpha,
\qquad
0\leq\widetilde h(y)\leq1-(1+\delta^2)^{-1/2}\leq\delta^2/2.
\]
We extend $\widetilde h$ by zero outside $D_{2R}^+$.

\smallskip
To estimate the BMO seminorm, let $y\in {D_{2R}^+}$, $x=\Psi(y)$, $0<\rho<\min\{R_0/2,R/4\}$, and set
\[
E=\Psi(B_\rho(y)\cap {D_{2R}^+}),\qquad S=B_{2\rho}(x)\cap {C_{2R}}.
\]
Since $\Psi$ and $\Phi$ have Lipschitz constants at most $2$,
\[
E\subset S\subset\Psi(B_{4\rho}(y)\cap {D_{2R}^+}),
\qquad \mu(S)\leq C(n,\alpha)\mu(E).
\]
The last inequality follows from the preceding bounds for $\omega$ and the doubling inequality in Lemma \ref{lm-doubling}-\textup{(i)} for $y_n^\alpha$ on ${D_{2R}^+}$. As $\det\nabla\Psi=1$, we obtain
\begin{align*}
&\fint_{B_\rho(y)\cap {D_{2R}^+}}
 |\lambda w-\overline{\lambda w}_{B_\rho(y)\cap {D_{2R}^+}}|\,d\omega
=\fint_E|\lambda u-\overline{\lambda u}_E|\,d\mu\\
&\quad\leq\frac{2}{\mu(E)}\int_E|\lambda u-\overline{\lambda u}_S|\,d\mu
\leq\frac{2\mu(S)}{\mu(E)}
\fint_S|\lambda u-\overline{\lambda u}_S|\,d\mu
\leq C(n,\alpha)M.
\end{align*}
Hence we take $M_*=C(n,\alpha)M$. Under $y=Rz$, the localization radius becomes
\[
\widehat R_0=\min\{R_0/(2R),1/4\}.
\]
{For $y\in D_R^+$ and $0<t<R\widehat R_0/2$, we have $D_t^+(y)\subset B_{2t}(y)\cap\bR^n_+$. Subtracting the weighted average of $\hat{\bA}$ on this ball in \eqref{theta-hat-A.def}, and using $\omega(x)\sim x_n^\alpha$ and Lemma \ref{lm-doubling}, we get
\[
\widehat\Theta_{t,y}(\hat{\bA},\omega)
\leq C(n,\alpha)\Theta_{2t,y}(\hat{\bA},\omega)
\leq 2C(n,\alpha)N_0\tau_2,
\]
where the last inequality follows from Lemma \ref{A-comparision}. Let $\overline\delta_0=\overline\delta_0(\Lambda_*,M_*,\alpha,p,n)$ be as in Theorem \ref{flat-thrm}. Choose $\tau_2\leq 1$ so that $2C(n,\alpha)N_0\tau_2<\overline\delta_0$ and $\tau_2^2/2\leq\overline\delta_0$.}}
{We apply Theorem \ref{flat-thrm} with localization radius $\widehat R_0/2$ to}
{\eqref{scaled-eqs} with $w_R(z)=w(Rz)/R$, $\lambda_R=R\lambda$, $\bG_R(z)=\bG(Rz)$, $\bA_R(z,s,\xi)=\hat{\bA}(Rz,s,\xi)$, and $\omega_R(z)=R^{-\alpha}\omega(Rz)$. Changing variables $y=Rz$, we obtain}
\[\|\nabla w (y)\|_{L^p(D_{R}^{+},\omega)}^{p} \leq C{\omega(D_{2R}^{+})^{1-p/2}}\|\nabla w(y)\|_{L^2(D_{2R}^{+},\omega)}^{p} + C{\|\bG\|_{L^p(D_{2R}^{+},\omega)}^{p}},\]
for $C={C(\Lambda,M,\alpha,p,n}, \hat{R}_0)>0$. Then, by applying the change of variables $x=\Psi(y)$ and using the definitions of $w$ and $\bG$ in Lemma \ref{change-of-variables-lm} together with \eqref{grads-psi-phi}, we infer \eqref{boundary-result}. This completes the proof.	 
\end{proof}

\subsection{Global estimates}\label{subsec:global}  We provide the proof of Theorem \ref{main-thrm}.
\begin{proof}[Proof of Theorem \ref{main-thrm}] Let $\bar{\rho}=\min\{\rho_0,r_0\}/32$ be as in the statement of the theorem, and let $\delta =\tfrac12\min\{\tau_1,\tau_2\}$ where $\tau_1$ is defined in Theorem \ref{interior-thrm} and $\tau_2$ is defined in Theorem \ref{boundary-thrm}. With the choice of $\delta$ and under the assumptions of Theorem \ref{main-thrm}, the conditions in \eqref{interior-eq} and \eqref{A-bmo-interior} are satisfied.  Since $\overline{\Omega}$ is compact, we choose finitely many boundary points $x_1, x_2, \ldots,x_N$ on $\partial \Omega$ such that
\[
\Omega\subset\Omega^{2\bar{\rho}}\cup\bigcup_{j=1}^N C_{3\rho_j}(x_j)
\] 
where $\rho_i \in [\bar{\rho}, \min\{\rho_0/24,r_0/6\})$ for all $i=1, 2, \ldots, N$. Here, note that $C_{\rho}(x)$ is defined in \ref{flatten-C}. We apply Theorem \ref{interior-thrm} to obtain
\[
\|\nabla u\|_{L^p(\Omega^{2\bar{\rho}},\mu)}^{p} \leq C\mu(\Omega)^{1-p/2}\|\nabla u\|_{L^2(\Omega,\mu)}^{p} + C \| \bF\|_{L^p(\Omega,\mu)}^{p}.\]
Similarly, for each $i\in\{1,\ldots,N\}$, applying Theorem \ref{boundary-thrm}, we obtain
	\[\|\nabla u\|_{L^p(C_{3\rho_i}(x_i),\mu)}^{p} \leq C\mu(\Omega)^{1-p/2}\|\nabla u\|_{L^2(\Omega,\mu)}^{p} + C \| \bF\|_{L^p(\Omega,\mu)}^{p}.\]
	Then, it follows that
\[\|\nabla u\|_{L^p(\Omega,\mu)}^{p} \leq {C\mu(\Omega)^{1-p/2}\|\nabla u\|}_{L^2(\Omega,\mu)}^{p} + C \| \bF\|_{L^p(\Omega,\mu)}^{p},\]
	where $C=C(\Lambda, M, p, \alpha, n, \bK, {R_0,\rho_0,r_0,\Omega,c_0,C_0})>0$ {also depends on the modulus of continuity of $\mu$ on ${\overline{\Omega^\rho}}$}.	
Moreover, by the energy estimate and H\"{o}lder's inequality, we get
\[\|\nabla u\|_{L^2(\Omega,\mu)} \leq {\Lambda}\|{\bF}\|_{L^2(\Omega,\mu)} \leq {\Lambda\mu(\Omega)^{1/2-1/p}}\|{\bF}\|_{L^p(\Omega,\mu)}.\]
Thus, combining this and the previous estimate, we obtain \[\|\nabla u\|_{L^p(\Omega,\mu)}^{p} \leq C \| \bF\|_{L^p(\Omega,\mu)}^{p}.\]
Applying the weighted Sobolev-Poincar\'e inequality (Proposition \ref{weighted-sobolev-poincare-lipschitz}) and the standard approximation argument, we infer \eqref{main-estimate2}. This completes the proof.

\end{proof}

\section*{Acknowledgement}
Preliminary results of this work were reported in Van Le's dissertation submitted to The University of Tennessee.  T. Phan's research is partially supported by Simons Foundation grant 769369,

\end{document}